\documentclass[12pt]{article}
\usepackage{amsmath}
\usepackage{graphicx,psfrag,epsf}
\usepackage{enumerate}
\usepackage{natbib}
\usepackage{url} 

\usepackage{mathrsfs,amsthm}
\usepackage{amssymb,multirow}
\usepackage{amsfonts}
\usepackage{dsfont}
\usepackage{fontenc}
\usepackage{graphicx}
\usepackage{subfigure}
\usepackage{bbm}
\usepackage{bm, mathabx}
\RequirePackage[OT1]{fontenc}

\newcommand{\indic}{\mathds{1}}

\newcommand{\X}{\mathds{X}}

\newcommand{\R}{\mathds{R}}
\newcommand{\E}{\mathds{E}}
\newcommand{\N}{\mathds{N}}
\renewcommand{\P}{\mathds{P}}

\newcommand{\Lcr}{\mathscr{L}}

\newcommand{\ddr}{\mathrm{d}}
\newcommand{\edr}{\mathrm{e}}

\def\simind{\stackrel{\mbox{\scriptsize{ind}}}{\sim}}
\def\simiid{\stackrel{\mbox{\scriptsize{iid}}}{\sim}}

\def\eqd{\stackrel{\mbox{\scriptsize{d}}}{=}}
\theoremstyle{plain}
\newtheorem{theorem}{\textsc{Theorem}}
\newtheorem{definition}{\textsc{Definition}}

\newtheorem{proposition}{\textsc{Proposition}}
\newtheorem{example}{\textsc{Example}}

\allowdisplaybreaks

\begin{document}

\def\spacingset#1{\renewcommand{\baselinestretch}%
	{#1}\small\normalsize} \spacingset{1}

{
	\title{\bf Latent nested nonparametric priors}
	\author{Federico Camerlenghi\\
		Department of Economics, Management and Statistics,\\ University of Milano--Bicocca\\
		and \\
		David B. Dunson\\
		Department of Statistical Science,\\ Duke University\\
		and\\
		Antonio Lijoi and Igor Pr\"unster\thanks{
		{A. Lijoi and I. Pr\"unster are supported by the European Research Council (ERC) through StG "N-BNP" 306406}}\hspace{.2cm}\\
		Department of Decision Sciences and BIDSA,\\ Bocconi University\\
		and \\
		Abel Rodr\'iguez \\
		Department of Applied Mathematics and Statistics, \\ University of California at Santa Cruz}
	\date{}
	\maketitle
}

\bigskip
\newpage

\begin{abstract}
	Discrete random structures are important tools in Bayesian nonparametrics and the resulting models have proven effective in density estimation, clustering, topic modeling and prediction, among others. In this paper, we consider nested processes and study the dependence structures they induce.  Dependence ranges between homogeneity, corresponding to full exchangeability, and maximum heterogeneity, corresponding to (unconditional) independence across samples. The popular nested Dirichlet process is shown to degenerate to the fully exchangeable case when there are ties across samples at the observed or latent level.  To overcome this drawback, inherent to nesting general discrete random measures, we introduce a novel class of latent nested  processes.  These are obtained by adding common and group-specific completely random measures and, then, normalising to yield dependent random probability measures.  We provide results on the partition distributions induced by latent nested processes, and develop an Markov Chain Monte Carlo sampler for Bayesian inferences. A test for distributional homogeneity across groups is obtained as a by product. The results and their inferential implications are showcased on synthetic and real data. 
\end{abstract}

\noindent%


\section{Introduction}\label{sec:intro}

Data that are generated from different (though related) studies, populations or experiments are typically characterised by some degree of heterogeneity. A number of Bayesian nonparametric models have been proposed to accommodate such data structures, but analytic complexity has limited understanding of the 
implied dependence structure across samples.  The spectrum of possible dependence ranges from homogeneity, corresponding to full exchangeability, to complete heterogeneity, corresponding to unconditional independence.  It is clearly desirable to construct a prior that can cover this full spectrum, leading to a posterior that can appropriately adapt to the true dependence structure in the available data.

This problem has been partly addressed in several papers. In \cite{lnp:2014} a class of random probability measures is defined in such a way that proximity to full exchangeability or independence is expressed in terms of a $[0,1]$--valued random variable. In the same spirit, a model decomposable into idiosyncratic and common components is devised in \cite{mqr:2004}. Alternatively, approaches based on P\'olya tree priors are developed in \cite{ma_wong:2011,hcgs:2015,filippi_olmo:2015}, while a multi--resolution scanning method is proposed in \cite{jacopo_li}. In \cite{arbishek_david} Dirichlet process mixtures are used to test homogeneity across groups of observations on a manifold.  A popular class of dependent nonparametric priors that fits this framework is the \textit{nested Dirichlet process} of \cite{ndp:2008}, which aims at clustering the probability distributions associated to $d$ populations. For $d=2$ this model is
\begin{equation}
  \label{eq:ndp}
    \begin{split}
    X_{i,1},X_{j,2}\,\mid\,(\tilde p_1,\tilde p_2)
    \: & \simind \: \tilde p_1\times\tilde p_2\\
    (\tilde p_1,\tilde p_2)\,\mid \,\tilde q\: &\sim \: \tilde q^2, 
    \qquad 
    \tilde q 
    =
    \sum_{i\ge 1}\omega_i\:\delta_{G_i}
 \end{split}
\end{equation}
where the random elements $\bm{X}_\ell=(X_{i,\ell})_{i\ge 1}$, for $\ell=1,2$, take values in a space $\X$, the sequences $(\omega_i)_{i\ge 1}$ and $(G_i)_{i\ge 1}$ are independent, with $\sum_{i\ge 1} \omega_i=1$ almost surely, and the $G_i$'s are i.i.d. random probability measures on $\X$ such that
\begin{equation}
  \label{eq:G_i}
  G_i=\sum_{t\ge 1}w_{t,i}\delta_{\theta_{t,i}},\qquad \theta_{t,i}\simiid P
\end{equation}
for some non--atomic probability measure $P$ on $\X$. In \cite{ndp:2008} it is assumed that $\tilde q$ and the $G_i$'s are realizations of Dirichlet processes while in \cite{abel_david_aoas} it is assumed they are from a generalised Dirichlet process introduced in  \cite{nils_gdp}. 
Due to discreteness of $\tilde q$, one has $\tilde p_1=\tilde p_2$ with positive probability allowing for clustering at the level of the populations' distributions and implying $\bm{X}_1\sim\bm{X}_2$ in such cases.  

The nested Dirichlet process has been widely used in a rich variety of applications, but it has an unappealing characteristic that provides motivation for this article.  In particular, if $\bm{X}_1$ and $\bm{X}_2$ share at least one value, then the posterior distribution of $(\tilde p_1,\tilde p_2)$ degenerates on $\{\tilde p_1=\tilde p_2\}$, forcing homogeneity across the two samples.  This occurs also in nested Dirichlet process mixture models in which the $X_{i,\ell}$ are latent, and is not specific to the Dirichlet process but is a consequence of nesting discrete random probabilities.  

To overcome this major limitation, we propose a more flexible class of {\em latent nested processes}, which preserve heterogeneity {\em a posteriori}, even when distinct values are shared by different samples.  Latent nested processes define $\tilde p_1$ and $\tilde p_2$ in \eqref{eq:ndp} as resulting from normalisation of an additive random measure model with common and idiosyncratic components, the latter with nested structure. Latent nested processes are shown to have appealing distributional properties. In particular, nesting corresponds, in terms of the induced partitions, to a convex combination of full exchangeability and unconditional independence, the two extreme cases.  This leads naturally to methodology for testing equality of distributions.


\section{Nested processes} \label{sec:gnp}
\subsection{Generalising nested Dirichlet processes via normalised random measures}
We first propose a class of nested processes that generalise nested Dirichlet processes by replacing the Dirichlet process components with a more flexible class of random measures.  The idea is to define $\tilde q$ in \eqref{eq:ndp} in terms of normalised \textit{completely random measures} on the space $\P$ of probability measures on $\X$.
Let $\tilde \mu$ be an almost surely finite completely random measure without fixed points of discontinuity, i.e. $\tilde\mu=\sum_{i\ge 1} J_i\:\delta_{G_i}$ where $G_i$ are i.i.d. random probability measures on $\X$ with some fixed distribution $Q$ on $\P$. The corresponding L\'evy measure on $\R^+\times\P$ is assumed to factorise as
\begin{equation}
\nu(\ddr s,\ddr p)=c\, \rho(s)\,\ddr s\: Q(\ddr p)
\label{eq:levy_1}
\end{equation}
where $\rho$ is some non--negative function such that $\int_0^\infty\min\{1,s\}\:\rho(s)\,\ddr s<\infty$ and $c>0$. Since such a $\nu$ characterises $\tilde\mu$ through its L\'evy-Khintchine representation
\begin{equation}
  \label{eq:expon_laplace}
  \E\Big[\edr^{-\lambda \tilde\mu(A)}\Big]=\exp\Big[-c\, Q(A)\,\int_{0}^\infty \Big(1-\edr^{-\lambda s}\Big)\,\rho(s)\,\ddr s \Big]=:\edr^{-c\, Q(A)\,\psi(\lambda)}
\end{equation}
for any measurable $A\subset \P$, we use the notation $\tilde\mu\sim\mbox{CRM}[\nu;\P]$. The function $\psi$ in \eqref{eq:expon_laplace} is also referred to as the \textit{Laplace exponent} of $\tilde\mu$. For a more extensive treatment of completely random measures, see \cite{Kingman:1993}. If one additionally assumes that $\int_0^\infty \rho(s)\:\ddr s=\infty$, then $\tilde \mu(\P)>0$ almost surely and we can define $\tilde q$ in \eqref{eq:ndp} as 
\begin{equation}
  \tilde q\,\eqd\,\frac{\tilde\mu}{\tilde\mu(\P)}
  \label{eq:nrmi_px}
\end{equation}
This is known as a \textit{normalised random measure with independent increments}, introduced in \cite{RLP:2003}, and is denoted as $\tilde q\sim\mbox{NRMI}[\nu ;\P]$. The baseline measure, $Q$, of $\tilde\mu$ in \eqref{eq:levy_1} is, in turn, the probability distribution of $\tilde q_0\sim\mbox{NRMI}[\nu_0;\X]$,  
with $\tilde q_0 \eqd \tilde\mu_0 / \tilde\mu_0(\X)$ and $\tilde \mu_0$ having L\'evy measure
\begin{equation}
\nu_0(\ddr s,\ddr x)=c_0\, \rho_0(s)\,\ddr s\: Q_0(\ddr x)
\label{eq:levy_2}
\end{equation}
for some non--negative function $\rho_0$ such that $\int_0^\infty\min\{1,s\}\,\rho_0(s)\,\ddr s<\infty$ and  $\int_0^\infty \rho_0(s)\:\ddr s=\infty$.
Moreover, $Q_0$ is a non--atomic probability measure on $\X$ and $\psi_0$ is the Laplace exponent of $\tilde\mu_0$. 
The resulting general class of nested processes is such that $(\tilde p_1,\tilde p_2)|\tilde q\sim \tilde q^2$ and is indicated by 
$(\tilde p_1,\tilde p_2)\sim
\mbox{NP}(\nu_0,\nu).$
The nested Dirichlet process of \cite{ndp:2008} is recovered by specifying $\tilde\mu$ and $\tilde\mu_0$ to be gamma processes, namely $\rho(s)=\rho_0(s)=s^{-1}\,\edr^{-s}$, so that both $\tilde q$ and $\tilde q_0$ are Dirichlet processes. 

\subsection{Clustering properties of nested processes}

A key property of nested processes is their ability to cluster both population distributions and data from each population.  In this subsection, we present results on: (i) the prior probability that $\tilde p_1=\tilde p_2$ and the resulting impact on ties at the observations' level; (ii) equations for mixed moments as convex combinations of fully exchangeable and unconditionally independent special cases; and (iii) a similar convexity result for partially exchangeable partition probability function. The probability distribution of an exchangeable partition depends only on the numbers of objects in each group; the exchangeable partition probability function is the probability of observing a particular partition as a function of the group counts. Partial exchangeability is exchangeability within samples; the partially exchangeable partition probability function depends only on the number of objects in each group that are idiosyncratic to a group and common. Simple forms for the partially exchangeable partition probability function not only provide key insights into the clustering properties but also greatly facilitate computation.

Before stating result (i), define
\[
\tau_q(u)=\int_0^\infty s^q\:\edr^{-us}\,\rho(s)\,\ddr s,\qquad
\tau_{q}^{(0)}(u)=\int_0^\infty s^q\:\edr^{-us}\,\rho_0(s)\,\ddr s,
\]
for any $u>0$, and agree that $\tau_0(u)\equiv\tau_0^{(0)} (u)\equiv 1$. 

\begin{proposition} \label{prp:coincidence} If $(\tilde p_1,\tilde p_2)\sim\mbox{\rm NP}(\nu_0,\nu)$, $c=Q(\P)$ and 
$c_0=Q_0(\X)$, then
  \begin{equation}
    \label{eq:prob_coincidence}
    \pi_1:=\P(\tilde p_1=\tilde p_2)=c\:\int_0^\infty u\:\edr^{-c\psi(u)}\:\tau_2(u)\:\ddr u
  \end{equation}
  and the probability that any two observations from the two samples coincide equals
  \begin{equation} \label{eq:ties}
  \P(X_{j,1}=X_{k,2})=\pi_1\:c_0\,\int_0^\infty u\:\edr^{-c_0\,\psi_0(u)} \:\tau_2^{(0)}(u)\,\ddr s>0.
  \end{equation} 
\end{proposition}
This result shows that the probability of $\tilde p_1$ and $\tilde p_2$ coinciding is positive, as desired, but also that this implies a positive probability of ties at the observations' level. Moreover, \eqref{eq:prob_coincidence} only depends on 
$\nu$ and  not 
{$\nu_0$}, since the latter acts on the $\X$ space.  In contrast, the probability that any two observations $X_{j,1}$ and $X_{k,2}$ from the two samples coincide given in \eqref{eq:ties} depends also on 
{$\nu_0$}.  If $(\tilde p_1,\tilde p_2)$ is a nested Dirichlet process, which corresponds to $\rho(s)=\rho_0(s)=\edr^{-s}/s$, one obtains $\pi_1=1/(c+1)$ and $\P(X_{1,1}=X_{1,2})=\pi_1/(c_0+1)$.

The following proposition [our result (ii)] provides a representation of mixed moments as a convex combination of full exchangeability and unconditional independence between samples.  
\begin{proposition} \label{prp:moment2} If {$(\tilde p_1,\tilde p_2)\sim\mbox{\rm NP}(\nu_0,\nu)$} and $\pi_1=\P (\tilde{p}_1=\tilde{p}_2)$ is as in \eqref{eq:prob_coincidence}, then 
\begin{equation}
\label{eq:2nd_moment}
\begin{split}
\E \Big[\int_{\P_\X^2} & f_1(p_1) f_2(p_2) \tilde{q}(\ddr p_1)
\tilde{q}(\ddr p_2)\Big] \\[4pt]
& = \pi_1 \int_{\P_{\X}} f_1 (p) f_2 (p) Q(\ddr p) + (1-\pi_1)
\int_{\P_{\X}} f_1(p) Q(\ddr p) \int_{\P_{\X}} f_2(p) Q(\ddr p)
\end{split}
\end{equation}
for all measurable functions $f_1,f_2 : \P \rightarrow \R^+$.
\end{proposition}
This convexity property is a key property of nested processes.

The component with weight $1-\pi_1$ in \eqref{eq:2nd_moment} accounts for heterogeneity among data from different populations and it is important to retain this component also \textit{a posteriori} in \eqref{eq:ndp}. Proposition~\ref{prp:moment2} is instrumental to obtain our main result (iii) characterizing the partially exchangeable random partition induced by $\bm{X}_1^{(n_1)}=(X_{1,1},\ldots,X_{n_1,1})$ and $\bm{X}_2^{(n_2)}=(X_{1,2},\ldots,X_{n_2,2})$ in \eqref{eq:ndp}. To fix ideas consider a partition of the $n_i$ data of sample $\bm{X}_i^{(n_i)}$ into $k_i$ specific groups and $k_0$ groups shared with sample $\bm{X}_j^{(n_j)}$ ($j\ne i$) with corresponding frequencies $\bm{n}_i=(n_{1,i},\ldots,n_{k_i,i})$ and $\bm{q}_i=(q_{1,i},\ldots,q_{k_0,i})$. For example, $\bm{X}_1^{(7)}=($0.5, 2, $-1$, 5, 5, 0.5, 0.5$)$ and $\bm{X}_2^{(4)}=($5, $-2$, 0.5, 0.5$)$ yield a partition of $n_1+n_2=11$ objects into $5$ groups of which $k_1=2$ and $k_2=1$ are specific to the first and the second sample, respectively, and $k_0=2$ are shared. Moreover, the frequencies are $\bm{n}_1=(1,1)$, $\bm{n}_2=(1)$, $\bm{q}_1=(3,2)$ and $\bm{q}_2=(2,1)$. Let us start by analyzing the two extreme cases. For the fully exchangeable case (in the sense of exchangeability holding true across both samples), one obtains the exchangeable partition probability function
\begin{equation}
  \label{eq:eppf_0}
  \begin{split}
    \Phi^{(N)}_{k}(\bm{n}_1,\bm{n}_2,\bm{q}_1+\bm{q}_2) &=
\frac{c_0^{k}}{\Gamma(N)} \int_0^{\infty} u^{N-1} \edr^{-c_0\psi_0(u)} \\
&\qquad\times\:
\prod_{j=1}^{k_1}\tau_{n_{j,1}}^{(0)} (u) \prod_{i=1}^{k_2}\tau_{n_{i,2}}^{(0)} (u)
\prod_{r=1}^{k_0}\tau_{q_{r,1}+q_{r,2}}^{(0)} (u)\: \ddr u
  \end{split}
\end{equation}
having set $N=n_1+n_2$, $k=k_0+k_1+k_2$ and  $|\bm{a}|=\sum_{i=1}^p a_i$ for any vector $\bm{a}=(a_1,\ldots,a_p)\in\R^p$ with $p\ge 2$.
The marginal exchangeable partition probability functions for the individual sample $\ell=1,2$ are
\begin{equation}
  \label{eq:margin_peppf}
    \Phi^{(n_\ell)}_{k_0+k_\ell}(\bm{n}_\ell,\bm{q}_\ell)
    =\frac{(c_0)^{k_0
        +k_\ell}}{\Gamma(n_\ell)} \:
 \int_0^{\infty} u^{n_\ell-1}\: \edr^{-c_0\,\psi_0(u)}
\prod_{j=1}^{k_\ell}\tau_{n_{j,\ell}}^{(0)} (u) \prod_{r=1}^{k_0}\tau_{q_{r,\ell}}^{(0)} (u) \:\ddr u
\end{equation}
Both \eqref{eq:eppf_0} and \eqref{eq:margin_peppf} hold true with the constraints $\sum_{j=1}^{k_\ell}n_{j,\ell}+\sum_{r=1}^{k_0}q_{r,\ell}=n_\ell$ and $1\le k_\ell+k_0\le n_\ell$, for each $\ell=1,2$. Finally, the convention $\tau_0 ^{(0)}\equiv 1$ implies that whenever an argument of the function $\Phi_k^{(n)}$ is zero, then it reduces to $\Phi_{k-1}^{(n)}$. For example, $\Phi_3^{(6)}(0,2,4)=\Phi_2^{(6)}(2,4)$.  Both \eqref{eq:eppf_0} and \eqref{eq:margin_peppf} solely depend on the L\'evy intensity of the completely random measure and can be made explicit for specific choices. We are now ready to state our main result (iii). 

\begin{theorem}\label{thm:EPPF_ndp}
	The random partition induced by the samples $\bm{X}_1^{(n_1)}$ and $\bm{X}_2^{(n_2)}$ drawn from $(\tilde p_1,\tilde p_2)\sim\mbox{\rm NP}(\nu_0,\nu)$, according to \eqref{eq:ndp}, is characterised by the partially exchangeable partition probability function
\begin{equation}
\label{eq:EPPF_nDP}
\begin{split}
\Pi^{(N)}_{k} & (\bm{n}_1,\bm{n}_2,\bm{q}_1,\bm{q}_2) =\pi_1\: \Phi^{(N)}_{k}(\bm{n}_1,\bm{n}_2,\bm{q}_1+\bm{q}_2)
\\[4pt]
& +  (1-\pi_1)\: \Phi^{(n_1 +|\bm{q}_1|)}_{k_0+k_1}(\bm{n}_1,\bm{q}_1) \Phi^{(n_2 + |\bm{q}_2|)}_{k_0+k_2}
(\bm{n}_2,\bm{q}_2) \indic_{\left\lbrace 0\right\rbrace}(k_0)
\end{split}
\end{equation}
\end{theorem}
The two independent exchangeable partition probability functions in the second summand on the right--hand side of \eqref{eq:EPPF_nDP} are crucial for accounting for the heterogeneity across samples. However, the result shows that one shared value, i.e. $k_0\ge 1$, forces the random partition to degenerate to the fully exchangeable case in \eqref{eq:eppf_0}.  Hence, a single tie forces the two samples to be homogeneous, representing a serious limitation of all nested processes including the nDP special case.  This result shows that degeneracy is a consequence of combining simple discrete random probabilities with nesting.  In the following section, we develop a generalisation that is able to preserve heterogeneity in presence of ties between the samples. 

\section{Latent nested processes} \label{sec:ncrm}

\begin{figure}[h!]
\begin{center}
\subfigure[]{\includegraphics[width=0.48\linewidth]{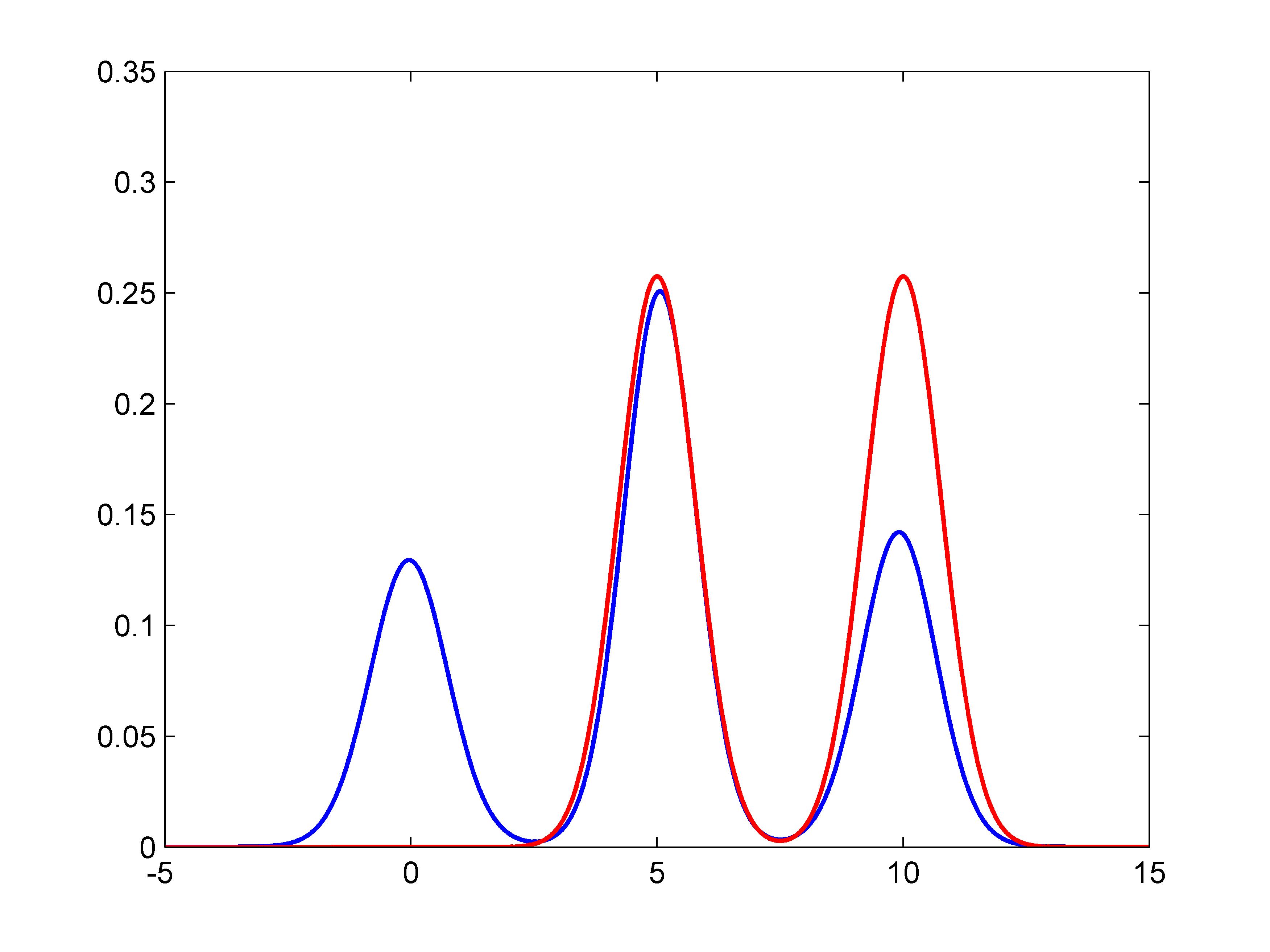}}
\subfigure[]{\includegraphics[width=0.48\linewidth]{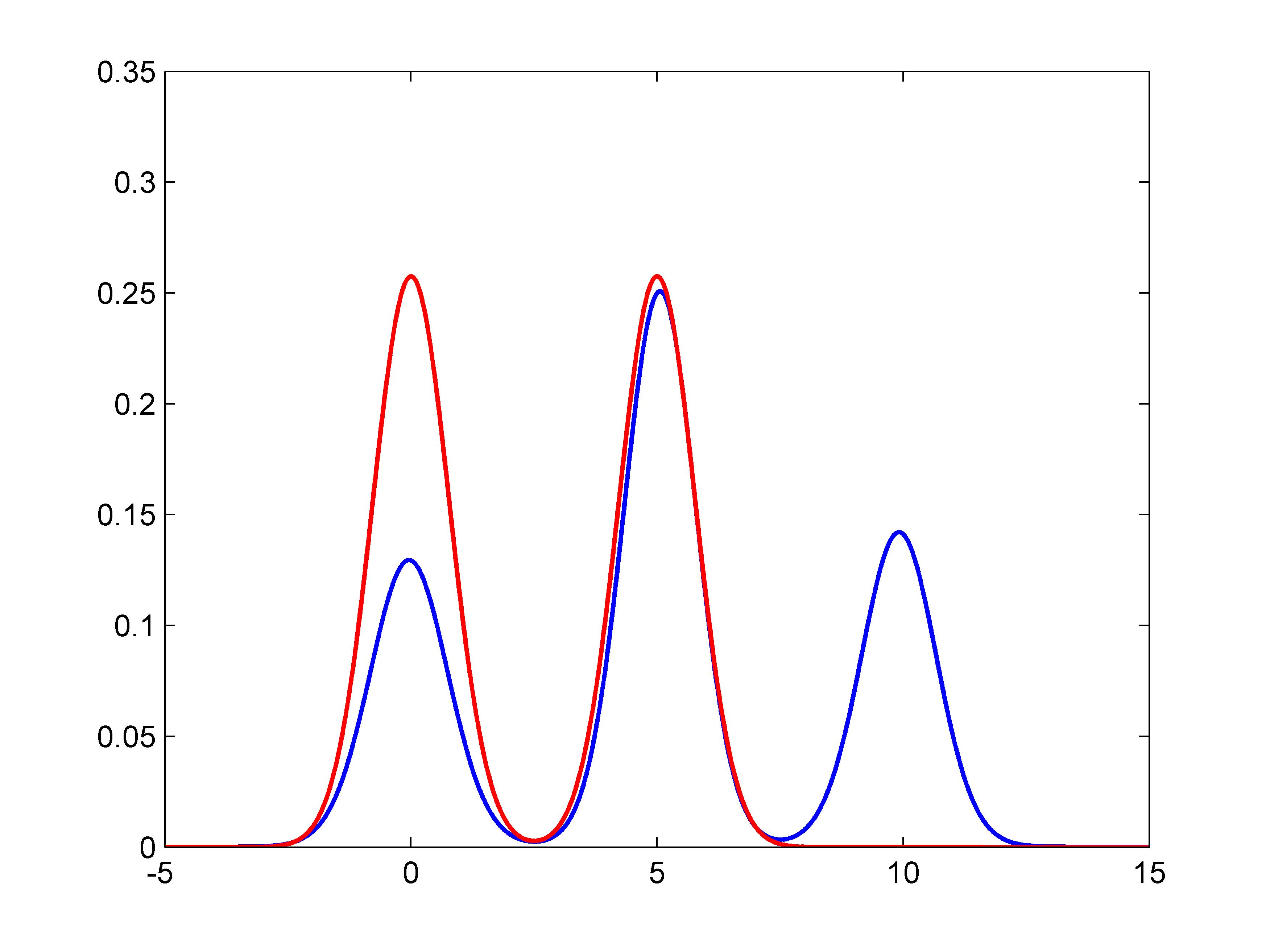}}
\end{center}
\begin{flushleft} 
\caption{{\rm Nested $\sigma$--stable mixture models: Estimated densities (blue) and true densities (red), for $\bm{X}_1^{(100)}$ in Panel (a) and for $\bm{X}_2^{(100)}$ in Panel (b).}}
\label{fig:nested_density_100}
\end{flushleft}
\end{figure}

To address degeneracy of the partially exchangeable partition probability function in \eqref{eq:EPPF_nDP}, we look for a model that, while still able to cluster random probabilities, can also take into account heterogeneity of the data in presence of ties between $\bm{X}_1^{(n_1)}$ and $\bm{X}_2^{(n_2)}$. The issue is relevant also in mixture models where $\tilde p_1$ and $\tilde p_2$ are used to model partially exchangeable latent variables such as, e.g., vectors of means and variances in normal mixture models. To see this, consider a simple density estimation problem, where two-sample data of sizes $n_1=n_2=100$ are generated from
\[
X_{i,1}\sim \frac{1}{2}\,\mbox{N}(5,\,\mbox{0.6})+
\frac{1}{2}\,\mbox{N}(10,\, \mbox{0.6})
\qquad X_{j,2}\sim \frac{1}{2}\,\mbox{N}(5,\, \mbox{0.6})+\frac{1}{2}\,\mbox{N}(0,\, \mbox{0.6}).
\]
This can be modeled by dependent normal mixtures with mean and variance specified in terms of a nested structure as in \eqref{eq:ndp}. The results, carried out by employing the algorithms detailed in Section \ref{sec:MCMC_algorithm}, show two possible outcomes:  either the model is able to estimate well the two bimodal marginal densities, while not identifying the presence of a common component, or it identifies the shared mixture component but does not yield a sensible estimate of the marginal densities, which both display three modes. The latter situation is displayed in Figure~\ref{fig:nested_density_100}: once the shared component $(5,\, \mbox{0.6})$ is detected, the two marginal distributions are considered identical as the whole dependence structure boils down to exchangeability across the two samples.

This critical issue can be tackled by a novel class of latent nested processes. Specifically, we introduce a model where the nesting structure is placed at the level of the underlying completely random measures, which leads to greater flexibility while preserving tractability.  In order to define the new process, let $\mathds{M}$ be the space of boundedly finite measures on $\X$ and $Q$ the probability measure on $\mathds{M}$ induced by $\tilde\mu_0\sim\mbox{CRM}[\nu_0;\X]$, where $\nu_0$ is as in \eqref{eq:levy_2}. Hence, for any measurable subset $A$ of $\X$ 
		\begin{equation*}
	\E\Big[\edr^{-\lambda\tilde\mu_0(A)}\Big]
	=\int_{\mathds{M}}\edr^{-\lambda\,
		m(A)}\, Q(\ddr m)=\exp \Big\{-
	c_0\,Q_0(A)\,
	\,\int_0^\infty\Big(1-\edr^{-\lambda s}\Big)\,\rho_0(s)\,\ddr s  \Big\}.
	\end{equation*}

\begin{definition}\label{def:nlp}
	{\rm Let $\tilde q\sim\mbox{NRMI}[\nu;\mathds{M}]$, with $\nu (\ddr s, \ddr m)= c \rho (s) \ddr s\, Q(\ddr m)$. Random probability measures $(\tilde p_1,\tilde p_2)$ are a \textit{latent nested process} if
	\begin{equation}
	\label{eq:def_nlp}
	\tilde p_\ell=\frac{\mu_{\ell}+\mu_S}{\mu_{\ell}(\X)+\mu_S(\X)}\qquad\ell=1,2,
	\end{equation}
	where $(\mu_1,\mu_2,\mu_S)\,|\,\tilde q\sim\tilde q^2\times \tilde q_S$ and $\tilde q_S$ is the law of a $\mbox{CRM}[\nu_0^*;\X]$, where $\nu_0^*=\gamma\,\nu_0$, for some $\gamma> 0$. Henceforth, we will use the notation $(\tilde p_1,\tilde p_2)\sim \mbox{LNP}(\gamma,\nu_0,\nu)$.}
\end{definition}
Furthermore, since
\begin{equation}
\tilde{p}_i= w_i \frac{\mu_i}{\mu_i(\X)} +(1-w_i) \frac{\mu_S}{\mu_S(\X)}, \quad
\text{where } \: w_i = \frac{\mu_i(\X)}{\mu_S(\X)+\mu_i(\X)},
\label{eq:gm_dependent}
\end{equation}
each $\tilde{p}_i$ is a mixture of two components: an idiosyncratic component $\mu_i/\mu_i(\X)$ and a shared component $\mu_S/\mu_S(\X)$. Here $\mu_S$ preserves heterogeneity across samples even when shared values are present. The parameter $\gamma$ in the intensity $\nu_0^*$ tunes the effect of such a shared CRM.  One recovers model \eqref{eq:ndp} as $\gamma\to 0$. 
A generalisation to nested completely random measures of the results given in Propositions \ref{prp:coincidence} and \ref{prp:moment2} is provided in the following proposition, whose proof is omitted.

\begin{proposition}\label{prp:crm_nested_coincidence}
If $(\mu_1,\mu_2)\,|\,\tilde q\sim \tilde q^2$, where $\tilde q\sim\mbox{\rm NRMI}[\nu;\mathds{M}]$ as in {\rm Definition~\ref{def:nlp}}, 
	then
\begin{equation}
  \label{eq:prob_coincid_crms}
  \pi_1^*=\P(\mu_1=\mu_2)
  =c\int_0^\infty u\:\edr^{-c\psi(u)}\:\tau_2(u)\:\ddr u
\end{equation}
and
\begin{equation}
  \label{eq:mixed_moment_crms}
  \begin{split}
    \E \Big[\int_{\mathds{M}^2} f_1(m_1)\, f_2(m_2)\: {\tilde q^2}(\ddr m_1,\ddr m_2) \Big]
    &=\pi_1^*\,\int_{\mathds{M}}f_1(m)\, f_2(m)\,Q(\ddr m)\\
    &\qquad\qquad +(1-\pi_1^*)\,\prod_{\ell=1}^2\int_{\mathds{M}}f_\ell(m)\,Q(\ddr m)
  \end{split}
\end{equation}
for all measurable functions $f_1,f_2:\mathds{M}\to\R^+$. 
\end{proposition}

\begin{proposition}\label{prp:probs_crms}
If 
$(\tilde p_1,\tilde p_2)\sim\mbox{\rm LNP}(\gamma,\nu_0,\nu)$, then $\P(\tilde p_1=\tilde p_2)=\P(\mu_1=\mu_2)$.
\end{proposition}

Proposition~\ref{prp:probs_crms}, combined with $\{\tilde p_1=\tilde p_1\}=\{\mu_1=\mu_2\}\cup(\{\tilde p_1=\tilde p_2\}\cap\{\mu_1\ne\mu_2\})$, entails $\P[\{\tilde p_1=\tilde p_2\}\cap\{\mu_1\ne \mu_2\}]=0$ namely
\[
\P(\{\tilde p_1=\tilde p_2\}\cap\{\mu_1 = \mu_2\})
+\P(\{\tilde p_1\ne \tilde p_2\}\cap\{\mu_1 \ne  \mu_2\})=1
\]  
and, then, the random variables $\indic\left\lbrace \tilde{p}_1 = \tilde{p}_2\right\rbrace$ and $\indic\left\lbrace  \mu_1= \mu_2 \right\rbrace$ coincide almost surely. As a consequence the posterior distribution of $\indic\left\lbrace  \mu_1= \mu_2 \right\rbrace$
can be readily employed to test equality between the distributions of the two samples. Further details are given in Section \ref{sec:illustrations}.

For analytic purposes, it is convenient to introduce an augmented version of the latent nested process, which includes latent 
indicator variables. In particular, $(X_{i,1},X_{j,2})\mid(\tilde p_1,\tilde p_2)\sim\tilde p_1\times\tilde p_2$, with $(\tilde p_1,\tilde p_2)\sim\mbox{LNP}(\gamma,\nu_0,\nu)$ if and only if 
\begin{equation}
\begin{split}
(X_{i,1 } , X_{j,2}) \mid (\zeta_{i,1} , \zeta_{j,2}, \mu_1, \mu_2,\mu_S)  
& \:\simind\: p_{\zeta_{1,i}} \times p_{2 \zeta_{2,j}} \\[4pt]
(\zeta_{i,1},\zeta_{j,2})\mid ( \mu_1, \mu_2 , \mu_S) 
& \: \sim\: { \rm Bern } (w_1) \times { \rm Bern } (w_2) \label{eq:lat_nested_CRM}\\[4pt]
( \mu_1, \mu_2 , \mu_S ) \mid (\tilde{q}, {\tilde{q}_S}) & \: \sim \: \tilde{q}^{2} \times {\tilde{q}_S}.
\end{split}
\end{equation}
The latent variables $\zeta_{i,\ell}$ indicate which random probability measure between $p_{\ell}$ and $p_0=p_S$  generates each $X_{i,\ell}$, for $i=1,\ldots,n_\ell$. 

\begin{theorem}\label{thm:peppf_nested_crms}
   The random partition induced by the samples $\bm{X}_1^{(n_1)}$ and $\bm{X}_2^{(n_2)}$ drawn from {$(\tilde p_1,\tilde p_2)\sim\mbox{\rm LNP}(\gamma,\nu_0,\nu)$}, 
   as in \eqref{eq:lat_nested_CRM}, is characterised by the partially exchangeable partition probability function
\begin{equation}
\label{eq:EPPF_nested_CRM2}
\begin{split}
\Pi_k^{(N)} & (\bm{n}_1 , \bm{n}_2 , \bm{q}_1 , \bm{q}_2) = \pi_1^* \frac{c_0^k (1+\gamma)^k}{
	\Gamma(N)}\\
&   \times\:\int_0^{\infty} s^{N-1}
e^{-(1+\gamma)c_0 \psi_0(s)} \prod_{\ell=1}^2 \prod_{j=1}^{k_\ell} \tau_{n_{j,\ell}}^{(0)} (s) 
\prod_{j=1}^{k_0} \tau_{q_{j,1}+q_{j,2}}^{(0)} (s) \ddr s  \\
& \qquad\qquad\qquad+
(1-\pi_1^*) \sum_{(\ast )}  I_2 (\bm{n}_1 , \bm{n}_2 , \bm{q}_1+ \bm{q}_2, \bm{\zeta}^*) 
\end{split}
\end{equation}
where 
\begin{align*}
&I_2 (\bm{n}_1 , \bm{n}_2 , \bm{q}_1+ \bm{q}_2, \bm{\zeta}^*) = \frac{c_0^k \gamma^{
k-\bar{k}
}}{\Gamma(n_1)
\Gamma (n_2)} \int_{0}^{\infty} \int_0^{\infty} u^{n_1-1}v^{n_2-1} e^{-\gamma c_0 \psi_0(u+v)-
c_0 (\psi_0(u)+ \psi_0(v))} \\
&\qquad\qquad\qquad \times \prod_{j=1}^{k_1} \tau_{n_{j,1}}^{(0)} (u+(1-\zeta_{j,1}^*)v)\prod_{j=1}^{k_2} \tau_{n_{j,2}}^{(0)} ((1-\zeta_{j,2}^*)u +v)
\prod_{j=1}^{k_0} \tau_{q_{j,1}+q_{j,2}}^{(0)} (u+v) \ddr u \ddr v
\end{align*}
and the sum in the second summand on the right hand side of \eqref{eq:EPPF_nested_CRM2} runs over all the possible labels $\bm{\zeta}^* \in \left\lbrace  0,1 \right\rbrace^{k_1+k_2}$.
\end{theorem}

The partially exchangeable partition probability function \eqref{eq:EPPF_nested_CRM2} is a convex linear combination of an exchangeable partition probability function corresponding to full exchangeability across samples and one corresponding to unconditional independence.  Heterogeneity across samples is preserved even in the presence of shared values.  The above result is stated in full generality, and hence may seem somewhat complex.  However, as the following examples show, when considering stable or gamma random measures, explicit expressions are obtained.  When $\gamma\to 0$ the expression \eqref{eq:EPPF_nested_CRM2} reduces to \eqref{eq:EPPF_nDP}, which means that the nested process is achieved as a special case. 

\begin{example}\label{exe:stable} 
Based on Theorem \ref{thm:peppf_nested_crms} we can derive an explicit expression of the partition structure of \textit{latent nested $\sigma$--stable processes}. Suppose $\rho(s)=\sigma\,s^{-1-\sigma}/\Gamma(1-\sigma)$ and $\rho_0(s)=\sigma_0\, s^{-1-\sigma_0}/\Gamma(1-\sigma_0)$, for some $\sigma$ and $\sigma_0$ in $(0,1)$. In such a situation it is easy to see that 
$\pi_1^*=1-\sigma$, $\tau_q^{(0)}(u)= \sigma_0(1-\sigma_0)_{q-1} u^{\sigma_0-q} $ and $\psi_0(u) = u^{\sigma_0}$. Moreover let $c_0=c=1$, since the total mass of a stable process 
is redundant under normalization. If we further set
\[
J_{\sigma_0, \gamma}(H_1,H_2;H):= \int_0^{1} \frac{w^{H_1-1} (1-w)^{
		H_2 -1}}{[\gamma +w^{\sigma_0}+(1-w)^{\sigma_0}]^H} \ddr w ,
\]
for any positive $H_1$, $H_2$ and $H$, and
\[
\xi_{a} (\bm{n}_1 , \bm{n}_2 , \bm{q}_1 + \bm{q}_2):=
\prod_{\ell=1}^2 \prod_{j=1}^{ k_\ell} 
(1-a)_{n_{j,\ell}-1} \prod_{j=1}^{k_0} (1-a)_{q_{j,1}+q_{j,2}-1} , 
\]
for any $a\in[0,1)$, then  
the partially exchangeable partition probability function in \eqref{eq:EPPF_nested_CRM2} may be rewritten as
\begin{multline*}
\Pi_k^{(N)} (\bm{n}_1 , \bm{n}_2 , \bm{q}_1 , \bm{q}_2) =
\sigma_0^{k-1} \Gamma(k) \xi_{\sigma_0} (\bm{n}_1 , \bm{n}_2 , \bm{q}_1 + \bm{q}_2) \left\{ \frac{(1-\sigma)}{\Gamma(N)} + \frac{\sigma}{\Gamma(n_1)\Gamma(n_2)} \right. \\ 
\left.
\times\:\sum_{(\ast )}  \gamma^{k-\bar{k}}\,
J_{\sigma_0, \gamma}(n_1-\bar{n}_1 +\bar{k}_1 \sigma_0 ,n_2-\bar{n}_2 +\bar{k}_2 \sigma_0;k)
\right\}.
\end{multline*}
The sum with respect to $\bm{\zeta}^*$ can be evaluated and it turns out that 
\begin{multline*}
\Pi_k^{(n)}(\bm{n}_1,\bm{n}_2,\bm{q}_1+\bm{q}_2)=
\frac{\sigma_0^{k-1} \Gamma(k)}{\Gamma(n)} \xi_{\sigma_0} (\bm{n}_1 , \bm{n}_2 , \bm{q}_1 + \bm{q}_2) 
\: \Big[1-\sigma + \sigma\gamma^{k_0}\: \frac{B(k_1\sigma_0,\, k_2\sigma_0)}{B(n_1,n_2)}\\
\times\: 
\int_0^1 \frac{\prod_{j=1}^{k_1}(1+\gamma w^{n_{j,1}-\sigma_0})\,
	\prod_{i=1}^{k_2}[1+\gamma(1-w)]^{n_{i,2}-\sigma_0}}
{\Big[\gamma+w^{\sigma_0}+(1-w)^{\sigma_0}\Big]^k}\:\mbox{Beta}(\ddr w; k_1\sigma_0,k_2\sigma_0)\,\Big]
\end{multline*}
where Beta$(\,\cdot\,;a,b)$ stands for the beta distribution with parameters $a$ and $b$, while $B(p,q)$ is the beta function with parameters $p$ and $q$. As it is well--known, $
\sigma_0^{k-1}\,\Gamma(k)\,\xi_{\sigma_0}(\bm{n}_1 , \bm{n}_2 , \bm{q}_1 + \bm{q}_2)/\Gamma(N)$
is the exchangeable partition probability function of a normalised $\sigma_0$--stable process. Details on the above derivation, as well as for the following example, can be found in the Appendix.
\end{example}

\begin{example}\label{exe:Dirichlet} 
Let $\rho(s)=\rho_0(s)=\edr^{-s}/s$. Recall that $\tau_q^{(0)}(u) = \Gamma(q)/(u+1)^q$ and $\psi_0(u)= \log(1+u)$, furthermore $\pi_1^*= 1/(1+c)$ by standard calculations.
From Theorem~\ref{thm:peppf_nested_crms} we obtain the partition structure of the \textit{latent nested Dirichlet process}
\begin{multline*}
\Pi_k^{(N)} (\bm{n}_1 , \bm{n}_2 , \bm{q}_1 , \bm{q}_2) = 
\xi_0(\bm{n}_1, \bm{n}_2, \bm{q}_1 +\bm{q}_2) c_0^k  
\left\{
\frac{1}{1+c} \,\frac{(1+\gamma)^k}{(c_0(1+\gamma))_N} \right. \\ 
+ \left.
\frac{c}{1+c} \:
\sum_{(\ast )} \frac{\gamma^{k-\bar{k}}
}{(\alpha)_{n_2} (\beta)_{n_1}
} {}_3F_2 (c_0+\bar{n}_2,\alpha, n_1; \alpha +n_2,
\beta+n_1;1 )   
\right\}
\end{multline*}
where $\alpha= (\gamma+1)c_0 +n_1-\bar{n}_1$, $\beta = c_0(2+\gamma)$ and ${}_3F_2$ is the generalised hypergeometric function. In the same spirit as in the previous example, the first element in the linear convex combination above $
c_0^k(1+\gamma)^k\:\xi_0(\bm{n}_1, \bm{n}_2, \bm{q}_1 +\bm{q}_2)/(c_0(1+\gamma))_N$ 
is nothing but the Ewens' sampling formula, i.e. the exchangeable partition probability function associated to the Dirichlet process whose base measure has total mass $c_0(1+\gamma)$.
\end{example}

\section{Markov Chain Monte Carlo algorithm} \label{sec:MCMC_algorithm}

We develop a class of Markov Chain Monte Carlo algorithms for posterior computation in latent nested process models relying on the partially exchangeable partition probability functions in 
Theorem~\ref{thm:peppf_nested_crms}, as 
they tended to be more effective. Moreover, the sampler is presented in the context of density estimation, where
\begin{align*}
  X_{j,\ell}\mid(\bm{\theta}_1^{(n_1)},\bm{\theta}_2^{(n_2)})
  \: &\simind\: h(\,\cdot\,;\theta_{j,\ell})\qquad \ell=1,2\\
  (X_{i,1},X_{j,2}) \mid(\bm{\theta}_1^{(n_1)},\bm{\theta}_2^{(n_2)})
  \: & \simind\: h(\,\cdot\,;\theta_{i,1})\,\times\, h(\,\cdot\,;\theta_{j,2})
\end{align*}
and the vectors $\bm{\theta}_\ell^{(n_\ell)}=(\theta_{1,\ell},\ldots,\theta_{n_\ell,\ell})$, for $\ell=1,2$ and with each $\theta_{i,\ell}$ taking values in $\Theta\subset\R^b$, are partially exchangeable and governed by a pair of $(\tilde p_1,\tilde p_2)$ as in \eqref{eq:lat_nested_CRM}. The discreteness of $\tilde p_1$ and $\tilde p_2$ entails ties among the latent variables $\bm{\theta}_1^{(n_1)}$ and $\bm{\theta}_2^{(n_2)}$ that give rise to $k=k_1+k_2+k_0$ distinct clusters identified by
\begin{itemize}
\item the $k_1$ distinct values specific to $\bm{\theta}_1^{(n_1)}$, i.e. not shared with $\bm{\theta}_2^{(n_2)}$. These are denoted as $\bm{\theta}_1^* :=(\theta_{1,1}^*, \dots, \theta_{k_1,1}^*)$, with corresponding frequencies $\bm{n}_1$ 
and labels $\bm{\zeta}^*_1$; 
\item the $k_2$ distinct values specific to $\bm{\theta}_2^{(n_2)}$, i.e. not shared with $\bm{\theta}_1^{(n_1)}$. These are denoted as 
$\bm{\theta}_2^* :=(\theta_{1,2}^*, \dots, \theta_{k_2,2}^*)$, with corresponding frequencies $\bm{n}_2$ 
and labels $\bm{\zeta}^*_2$; 
\item the $k_0$ distinct values shared by $\bm{\theta}_1^{(n_1)}$ and $\bm{\theta}_2^{(n_2)}$. These are denoted as 
$\bm{\theta}_0^* :=(\theta_{1,0}^*, \dots, \theta_{k_0,0}^*)$, with $\bm{q}_\ell$ 
being their frequencies in $\bm{\theta}_\ell^{(n_\ell)}$
and shared labels $\bm{\zeta}^*_0$. 
\end{itemize}
As a straightforward consequence of Theorem~\ref{thm:peppf_nested_crms}, one can determine the joint distribution of the data $\bm{X}$, the corresponding latent variables $\bm{\theta}$ and labels $\bm{\zeta}$ as follows
\begin{equation}\label{eq:joint}
  f(\bm{x}\mid\bm{\theta})\,\Pi_k^{(N)}(\bm{n}_1,\bm{n}_2,\bm{q}_1,\bm{q}_2)\,\prod_{\ell=0}^2\prod_{j=1}^{k_\ell}Q_0(\ddr\theta_{j,\ell}^*)
\end{equation}
where $\Pi_k^{(N)}$ is as in \eqref{eq:EPPF_nested_CRM2} and, for $C_{j,\ell}:=\{i:\: \theta_{i,\ell}=\theta_{j,\ell}^*\}$ and $C_{r,\ell,0}:=\{i:\: \theta_{i,\ell}=\theta_{r,0}^*\}$,
\[
f(\bm{x}\mid\bm{\theta})=\prod_{\ell=1}^2 \prod_{j=1}^{k_\ell}\:\prod_{i\in C_{j,\ell}}h(x_{i,\ell};\theta_{j,\ell}^*)\,\prod_{r=1}^{k_0}
\:\prod_{i\in C_{r,\ell,0}}\,h(x_{i,\ell};\theta_{r,0}^*).
\]
We do now specialise \eqref{eq:joint} to the case of latent nested $\sigma$--stable processes described in Example~\ref{exe:stable}. The Gibbs sampler is described just for sampling $\bm{\theta}_1^{(n_1)}$, since the structure is replicated for $\bm{\theta}_2^{(n_2)}$. To simplify the notation, $v^{-j}$ denotes the random variable $v$ after the removal of $\theta_{j,1}$. Moreover, with $\bm{T}=(\bm{X},\bm{\theta},\bm{\zeta},\sigma,\sigma_0,\bm{\phi})$, we let $\bm{T}_{-\theta_{j,1}}$ stand for $\bm{T}$ after deleting $\theta_{j,1}$, $I=\indic\{\tilde p_1=\tilde p_2\}$ and $Q_j^*(\ddr\theta)=h(x_{j,1};\theta)\, Q_0(\ddr\theta)/\int_\Theta h(x_{j,1};\theta)\, Q_0(\ddr\theta)$. Here $\bm{\phi}$ denotes a vector of hyperparameters entering the definition of the base measure $Q_0$. The updating structure of the Gibbs sampler is as follows
\medskip

\noindent (1) Sample $\theta_{j,1}$ from 
\begin{align*}
\P & (\theta_{j,1}\in \ddr \theta\, |\, \bm{T}_{-\theta_{j,1}}, I=1)= w_0 Q_{j,1}^* (\ddr \theta) +
\sum_{\left\lbrace i : \, \zeta_{i,0}^{*,-j}=\zeta_{j,1}\right\rbrace} w_{i,0} \delta_{\left\lbrace \theta_{i,0}^{*,-j}\right\rbrace}  (\ddr \theta)\\
&\quad +
\sum_{\left\lbrace i : \, \zeta_{i,1}^{*,-j}=\zeta_{j,1}\right\rbrace} w_{i,1} \delta_{\left\lbrace \theta_{i,1}^{*,-j}\right\rbrace} (\ddr \theta)+
\sum_{\left\lbrace i : \, \zeta_{i,2}^{*,-j}=\zeta_{j,1}\right\rbrace} w_{i,2} \delta_{\left\lbrace \theta_{i,2}^{*,-j}\right\rbrace} (\ddr \theta) \\[4PT]
\P & (\theta_{j,1}\in \ddr \theta \, |\, \bm{T}_{-\theta_{j,1}}, I=0)= w_0' Q_{j,1}^* (\ddr \theta)+
\sum_{\left\lbrace i : \, \zeta_{i,1}^{*,-j}=\zeta_{j,1}\right\rbrace} w_{i,1}' \delta_{\left\lbrace \theta_{i,1}^{*,-j}\right\rbrace} (\ddr \theta) \\
& \quad +
\indic_{\left\lbrace 0\right\rbrace}(\zeta_{j,1})\Big[\sum_{\left\lbrace i : \, \zeta_{i,2}^{*,-j}=0\right\rbrace} w_{i,2}' \delta_{\left\lbrace \theta_{i,2}^{*,-j}\right\rbrace} (\ddr \theta) +
\sum_{r=1}^{k_0} w_{r,0}' \delta_{\left\lbrace \theta_{r,0}^{*,-j}\right\rbrace}  (\ddr \theta)\Big]
\end{align*}
where 
\[
w_0 \:\propto\:\frac{\gamma^{1-\zeta_{j,1}}\sigma_0\,k^{-r} }{1+\gamma} h(x_{j,1} ; \theta),
\quad
w_{i,\ell}\: \propto\: (n_{i,\ell}^{-j}-\sigma_0) h(x_{j,1}; \theta_{i,\ell}^{*,-j}) \quad \ell=1,2
\]
\[
w_{i,0} \: \propto\: (q_{i,1}^{-j}+q_{i,2}^{-j}-\sigma_0) h(x_{j,1} ; \theta_{i,0}^{*,-j})
\]
and, with $a_1=n_1-(\bar{n}_1^{-j} +\zeta_{j,1}) +\bar{k}_1^{-j}\sigma_0$ and $a_2=n_2 -\bar{n}_2 +\bar{k}_2 \sigma_0$, one further has
\begin{align*}
w_0' & \:\propto\: \gamma^{1-\zeta_{j,1}}\sigma_0
k^{-j} J_{\sigma_0}(a_1+\zeta_{j,1}\sigma_0 , a_2 ; k^{-j} +1)\, h(x_{j,1} ; \theta), \\[4pt]
w_{i,\ell}' & \:\propto \:
J_{\sigma_0}(a_1, a_2 ; k^{-j} )\,
(n_{i,\ell}^{-j}-\sigma_0)\, h(x_{j,\ell}; \theta_{j,\ell}^{*,-j}) \qquad\qquad\qquad \ell=1,2,  \\[4pt]
w_{i,0}' & \:\propto\: J_{\sigma_0}(a_1, a_2 ; k^{-j} )\,
 (q_{i,1}^{-j}+q_{i,2}^{-j}-\sigma_0)\, h(x_{j,1} ; \theta_{i,0}^{*,-j}).
\end{align*}
\smallskip

\noindent (2) Sample $\zeta_{j,1}^*$ from
\begin{align*}
  \P(\zeta_{j,1}^*=x\mid \bm{T}_{-\zeta_{j,1}^*}, I=1 )
  &= \frac{\gamma^{1-x}}{1+\gamma}\\[4pt]
  \P (\zeta_{j,1}^*=x\mid \bm{T}_{-\zeta_{j,1}^*}, I=0 )
  &\:\propto \: \gamma^{k-k_x-\bar{k}_0-\bar{k}_2}  J_{\sigma_0}(n_1 -n_x +k_x\sigma_0,n_2-\bar{n}_2 +\bar{k}_2 \sigma_0;k)
\end{align*}
where $x\in\{0,1\}$, $k_x:= x+|\bm{\zeta}_{1}^{*,-j}| $
and $n_x=n_{j,1}x +|\bm{\zeta}_{1}^{*,-j}\odot\bm{n}_1^{-j}|$, where $\bm{a} \odot \bm{b}$ denotes the component--wise product between two vectors $\bm{a}, \bm{b}$. Moreover, it should be stressed that, conditional on $I=0$, the labels $\zeta_{r,0}^*$ are degenerate at $x=0$ for each $r=1,\ldots,k_0$. 
\medskip

\noindent (3) Update $I$ from 
\[
  \P(I=1\mid \bm{T})=1-\P(I=0\,|\,\bm{T})
  = \frac{ (1-\sigma)B(n_1,n_2)}{(1-\sigma)B(n_1,n_2) +\sigma J_{\sigma_0} (\bar{a}_1,\bar{a}_2; k) (1+\gamma)^k}
\]
where $\bar{a}_1=n_1 -\bar{n}_1 +\bar{k}_1\sigma_0$ and $\bar{a}_2=n_2-\bar{n}_2 +\bar{k}_2 \sigma_0$. This sampling distribution holds true whenever $\bm{\theta}_1^{(n_1)}$ and $\bm{\theta}_2^{(n_2)}$ do not share any value $\theta_{j,0}^*$ with label $\zeta_{j,0}^*=1$. If this situation occurs, then $\P(I=1\mid\bm{T})=1$. 
\medskip

\noindent (4) Update $\sigma$ and $\sigma_0$ from
\begin{align*}
  f(\sigma_0\,|\,\bm{T}_{-\sigma_0},I)
  & \:\propto\: J_{\sigma_0}^{1-I}(\bar{a}_1 , \bar{a}_2 ; k) \:\sigma_0^{k-1} {\kappa_0 
  	(\sigma_0)}  \prod_{\ell=1}^2 \prod_{j=1}^{k_\ell}(1-\sigma_0)_{n_{j,\ell}-1} \prod_{r=1}^{k_0} 
   (1-\sigma_0)_{q_{r,1}+q_{r,2} -1}\\[4pt]
f(\sigma \,|\, \bm{T}_{-\sigma}, I)
    & \:\propto\: \kappa (\sigma) \left[ (1-\sigma) \indic_{\{ 1 \}} (I) +
   \sigma \indic_{\{ 0 \}} (I) \right]
\end{align*}
where $\kappa$ and $\kappa_0$ are the priors for $\sigma$ and $\sigma_0$, respectively. 
\medskip

\noindent (5) Update $\gamma$ from
\begin{equation*}
   f(\gamma \,|\, \bm{T}_{-\gamma}, I)
   \:\propto\: \gamma^{k-\bar{k}} \, 
  {g(\gamma)}\:\Big[\frac{1-\sigma}{(1+\gamma)^k}\: \indic_{\left\lbrace  1 \right\rbrace} (I) + \sigma \, J_{\sigma_0} (\bar{a}_1, \bar{a}_2; k)\,\indic_{\{ 0 \}} (I) \Big]
\end{equation*}
where $g$ is the prior distribution for $\gamma$.
\medskip

Finally, the updating of the hyperparameters depends on the specification of $Q_0$ that is adopted. They will be displayed in the next section, under the assumption that $Q_0$ is a normal/inverse--Gamma.

The evaluation of the integral $J_{\sigma_0} (h_1, h_2;h)$ is essential for the implementation of the Markov Chain Monte Carlo procedure. This can be accomplished through numerical methods based on quadrature. However, computational issues arise when  $h_1$ and $h_2$ are both less than $1$ and the integrand defining $J_{\sigma_0}$ is no longer bounded, although still integrable. For this reason we propose a plain Monte Carlo approximation of $J_{\sigma_0}$ based on observing that
\[
J_{\sigma_0} (h_1, h_2;h)= B(h_1,h_2) \: \E \Big\{ \frac{1}{[\gamma +W^{\sigma_0}+ (1-W)^{\sigma_0}]^h} \Big\},
\]
with $W \sim {\rm Beta} (h_1,h_2)$. Then generating an i.i.d. sample $\left\lbrace W_i \right\rbrace_{i=1}^{L}$ of length $L$, with $W_i \sim W$, 
we get the following approximation
\[
J_{\sigma_0} (h_1, h_2;h) \approx B(h_1,h_2) \frac{1}{L}\sum_{i}^L\frac{1}{[\gamma +W_i^{\sigma_0}+ (1-W_i)^{\sigma_0}]^h} .
\]

\section{Illustrations}\label{sec:illustrations}

The algorithm introduced in Section \ref{sec:MCMC_algorithm} is employed here to estimate dependent random densities. 
Before implementation, we need first to complete the model specification of our latent nested model \eqref{eq:def_nlp}.  Let $\Theta= \R \times \R^+$
and  $h(\cdot; (M,V))$ be Gaussian with mean $M$ and variance $V$.
Moreover, as customary, 
$Q_0$ is assumed to be a normal/inverse--Gamma distribution
\[
Q_0(\ddr M, \ddr V) = Q_{0,1} (\ddr V) Q_{0,2} (\ddr M | V)
\]
with $Q_{0,1} $ an inverse--Gamma probability distribution with parameters $(s_0, S_0) $ 
and $Q_{0,2} $ a Gaussian with mean $m$ and variance $\tau V$. Furthermore, the hyperpriors are 
\begin{equation*}
\tau^{-1} \sim {\rm Gam}(w/2, W/2), \qquad
m  \sim {\rm N} (a,A),
\end{equation*}
for some real parameters $w>0, W>0, A>0$ and $a\in \R$. In the simulation  studies we have set
$(w,W)= (1,100)$, $(a,A)=((n_1 \bar{X} + n_2 \bar{Y})/(n_1+n_2),2)$.
The parameters $ \tau$ and $m$ are updated on the basis of their full conditional distributions, which can be easily derived, and correspond to
\begin{align*}
\Lcr (\tau | \bm{T}_{-\tau}, I) & \sim {\rm IG} \Big( \frac{w}{2}+ \frac{k}{2} ,
\frac{W}{2}+ \sum_{i=0}^2 \sum_{j=1}^{k_i} \frac{(M_{i,j}^* -m)^2}{2 V_{i,j}^*}\Big),\\
\Lcr (m | \bm{T}_{-m}, I) & \sim {\rm N} \Big( \frac{R}{D} , \frac{1}{D} \Big)
\end{align*}
where
\[
R= \frac{a}{A}+\sum_{i=0}^2\sum_{j=1}^{k_i} \frac{M_{i,j}^*}{\tau V_{i,j}^*}, 
\quad 
D= \frac{1}{A}+ \sum_{i=0}^2 \sum_{j=1}^{k_i}  \frac{1}{\tau V_{i,j}^*}.
\]
The model specification is completed by choosing uniform prior distributions for $\sigma_0$ and $\sigma$. In order to overcome the possible slow mixing of the P\'olya urn sampler, we include the acceleration step of \cite{macea_94} and \cite{westetal_94}, which consists in resampling the distinct values
$ (\theta_{i,j}^*)_{j=1}^{k_i} $, for $i=0,1,2$, at the end of every iteration.
The numerical outcomes displayed in the sequel are based on
$50,000$ iterations after $50,000$ burn--in sweeps. 

Throughout we assume the data $\bm{X}_1^{(n_1)}$ and $\bm{X}_2^{(n_2)}$ to be independently  generated by two densities $f_1$ and $f_2$.  These will be estimated jointly through the MCMC procedure and the borrowing of strength phenomenon should then allow improved performance.
An interesting byproduct of our analysis is the possibility to examine the 
clustering structure of each distribution, namely the number of components of each mixture.
Since the expression of the pEPPF \eqref{eq:EPPF_nested_CRM2} consists of two terms, in order to carry out 
posterior inference we have defined the random variable $I= \mathds{1}_{\{\mu_1=\mu_2\}}$.
This random variable allows to test whether the two samples come from the same distribution or not,
since $I= \indic_{\{ \tilde{p}_1 = \tilde{p}_2 \}}$ almost surely
(see also Proposition \ref{prp:probs_crms}).
Indeed, if interest lies in testing
\[
H_0 : \: \tilde{p}_1=\tilde{p}_2 \quad \text{versus} \quad H_1 :  \: \tilde{p}_1 \not =\tilde{p}_2,
\]
based on the Markov Chain Monte Carlo output, it is straightforward to compute an approximation of the Bayes factor
\[
\mathrm{BF} =   \frac{\P(\tilde{p}_1=\tilde{p}_2 | \bm{X})} 
	{\P(\tilde{p}_1\not=\tilde{p}_2 | \bm{X} )}
	\: \frac{\P(\tilde{p}_1\not = \tilde{p}_2)}{\P (\tilde{p}_1 = \tilde{p}_2)} 
	= \frac{\P(I=1| \bm{X})}
	{\P(I=0| \bm{X})} 
	\: \frac{\P (I=0)}{\P(I=1)}  
\]
leading to acceptance of the null hypothesis if $\mathrm{BF}$ is sufficiently large. In the following we first consider simulated datasets generated from normal mixtures and then we analyse the popular Iris dataset.

\subsection{Synthetic examples}\label{subsec:synthetic_data}

We consider three different simulated scenarios, 
where $\bm{X}_1^{(n_1)}$ and $\bm{X}_2^{(n_2)}$
are independent and identically distributed draws from densities that are both two component mixtures of normals. In both cases $(s_0, S_0)=(1,1)$ and the sample size is $n=n_1=n_2=100$.

First consider a scenario where $\bm{X}_1^{(n_1)}$ and $\bm{X}_2^{(n_2)}$ are drawn from the same density
\[
X_{i,1} \sim X_{j,2} \sim \frac{1}{2}\, {\rm N}(0,1)+\frac{1}{2}\, {\rm N}(5,1). 
\]
The posterior distributions for the number of mixture components, respectively denoted by $K_1$ and $K_2$ for the two samples, and for the number of shared components, denoted by $K_{12}$, are reported in Table \ref{tab:simulated_posterior}. The maximum a posteriori estimate is highlighted in bold. The model is able to detect the correct number of components for each distribution as well as the correct number of components shared across the two mixtures. The density estimates, not reported here, are close to the true data generating densities. The Bayes factor to test equality between the distributions of $\bm{X}_1^{(n_1)}$ and $\bm{X}_2^{(n_2)}$ has been approximated through the Markov Chain Monte Carlo output and coincides with $\mathrm{BF} = \mbox{5.85}$, providing evidence in favor of the null hypothesis.  

\begin{table}[h!]
	\begin{center} 
			\begin{tabular}{|c|l|c|c|c|c|c|c|c|c|}\hline 
				scen. & $\#$ comp.  & 0 & 1 & 2 & 3 & 4 & 5 & 6 & $\ge 7$ \\\hline \hline
				\multirow{3}{*}{I} & $K_1$ 
				& 0 & 0  & \bf{0.638}  & 0.232  & 0.079  & 0.029 & 0.012 & 0.008     \\   	
				& $K_2$ 
				& 0 &  0 &  \bf{0.635}  & 0.235  & 0.083 & 0.029 & 0.011 & 0.007      \\
				& $K_{12}$ 
				&  0  & 0  & \bf{0.754}  & 0.187   & 0.045 & 0.012 & 0.002 & 0.001     \\\hline
				\multirow{3}{*}{II} & $K_1$ 
				& 0 & 0 & \bf{0.679} & 0.232 & 0.065 & 0.018 &  0.004 & 0.002       \\   	
				& $K_2$ 
				& 0 & 0 & \bf{0.778} & 0.185 & 0.032 & 0.004 & 0.001 & 0        \\
				& $K_{12}$ 
				& 0 & \bf{0.965} & 0.034 & 0.001 & 0 & 0 & 0 & 0 \\\hline
				\multirow{3}{*}{III} & $K_1$  
				& 0 & 0 & \bf{0.328} & 0.322 & 0.188 & 0.089 &  0.041 & 0.032      \\   	
				& $K_2$ 
				& 0 & 0 & \bf{0.409} & 0.305 & 0.152 & 0.073 & 0.034 & 0.027        \\
				& $K_{12}$ 
				& 0 & 0.183 & \bf{0.645} & 0.138 & 0.027 & 0.006 & 0.001 & 0 \\\hline
			\end{tabular}
		\begin{flushleft}
			\caption{Simulation study: Posterior distributions of the number of components in the first sample ($K_1$), in the second sample ($K_2$) and shared by the two samples ($K_{12}$) corresponding to the three scenarios. The posterior probabilities corresponding to the MAP estimates are displayed in bold.}
			\label{tab:simulated_posterior}
		\end{flushleft}
	\end{center}
\end{table}

\begin{figure}[h!]
	\begin{center}
		\subfigure[]{\includegraphics[width=0.48\linewidth]{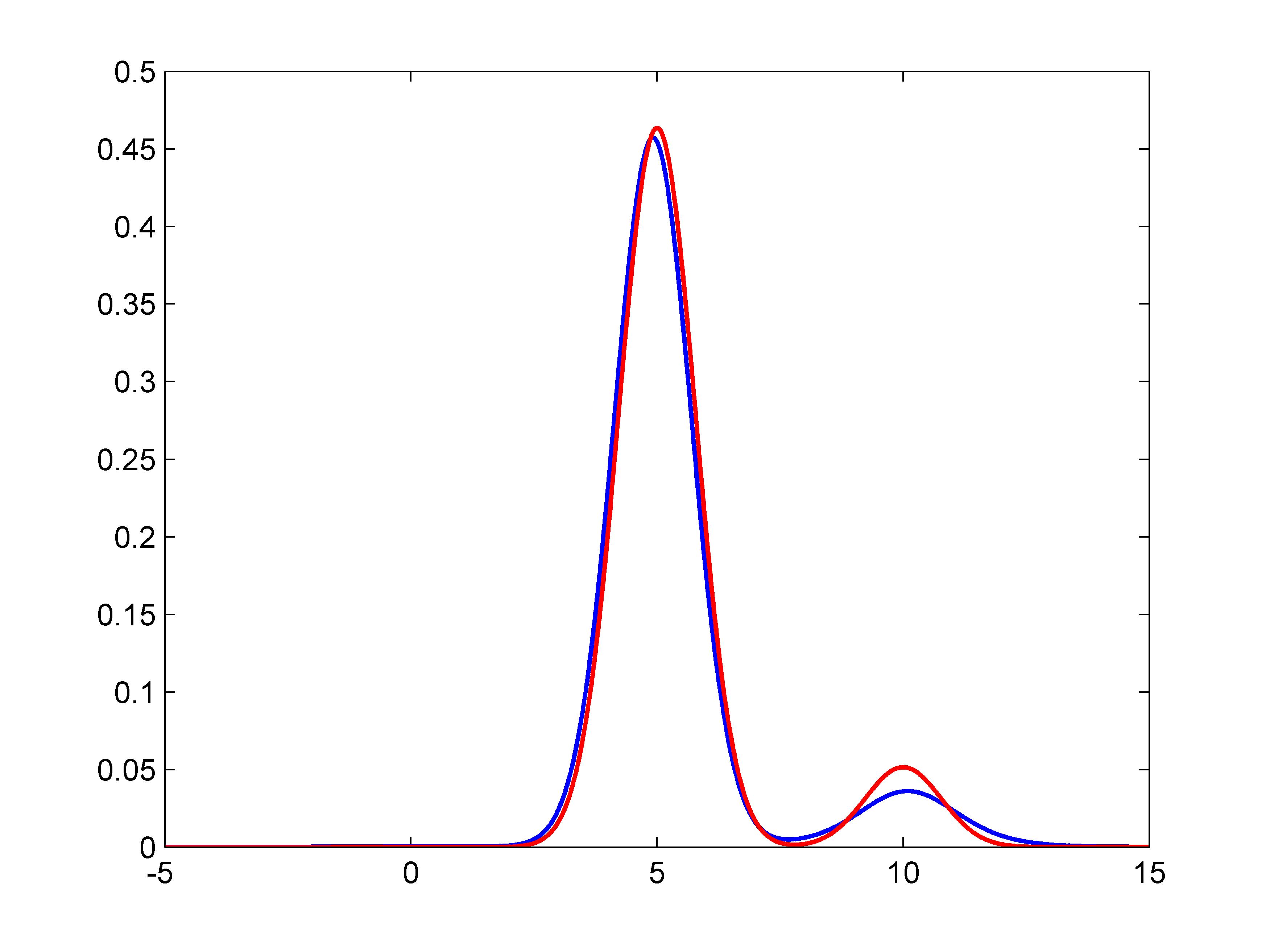} \label{fig: marginali_peso09_density_100a}}
		\subfigure[]{\includegraphics[width=0.48\linewidth]{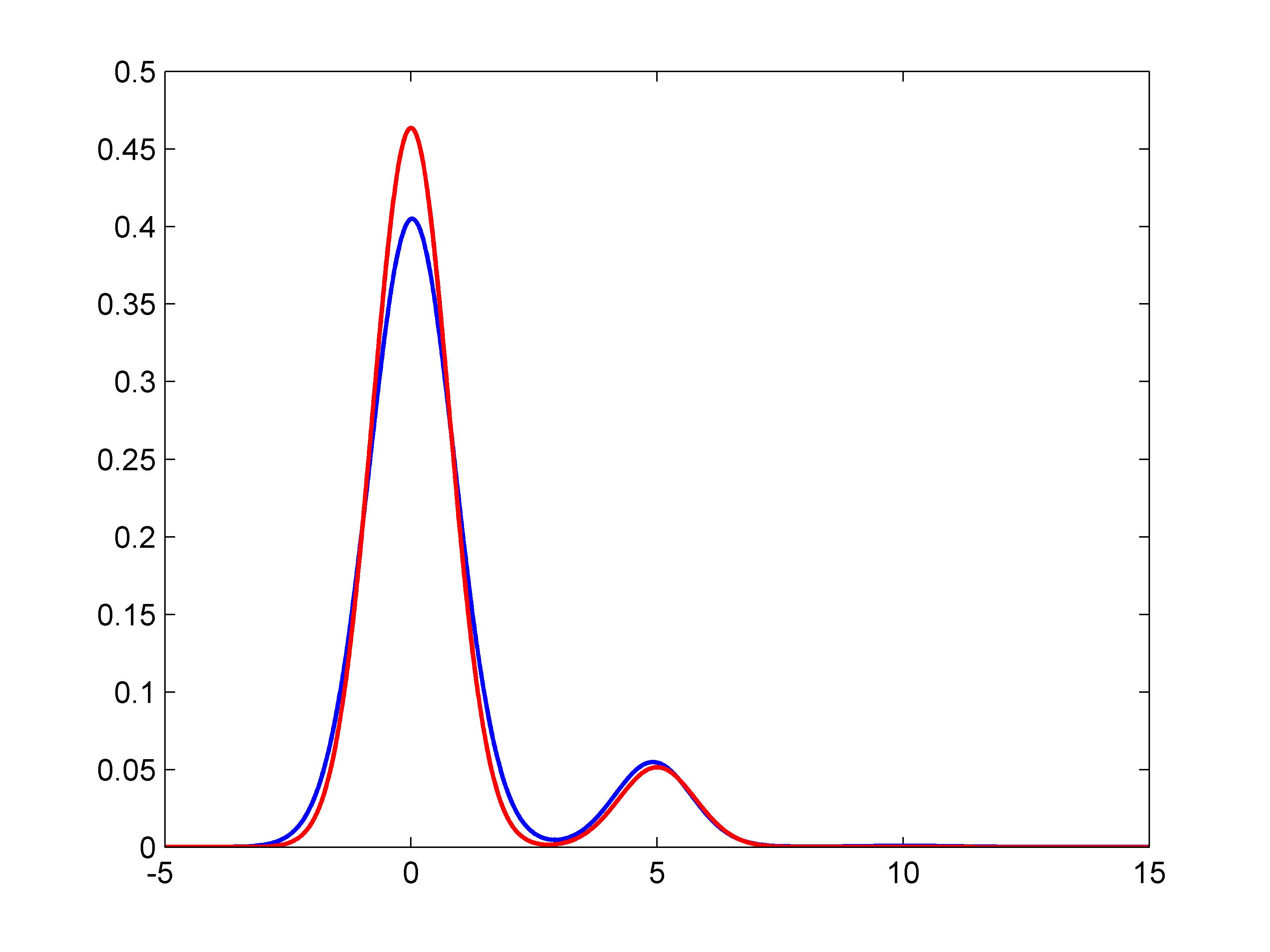} \label{fig: marginali_peso09_density_100b}}
	\end{center}
	\begin{flushleft} 
		\caption{Estimated densities (blue) and true densities (red) for $\mathbf{X}$ in Panel (a) and $\mathbf{Y}$ in Panel (b).}
		\label{fig: marginali_peso09_density_100}
	\end{flushleft}
\end{figure}

Scenario II corresponds to samples $\bm{X}_1^{(n_1)}$ and $\bm{X}_2^{(n_2)}$ generated, respectively, from
\[
X_{i,1} \sim \mbox{0.9}\, {\rm N}(5,\mbox{0.6})+\mbox{0.1}\, {\rm N}(10,\mbox{0.6}) \quad
X_{j,2} \sim \mbox{0.1}\, {\rm N}(5,\mbox{0.6})+\mbox{0.9}\, {\rm N}(0,\mbox{0.6}).
\]
Both densities have two components but only one in common, i.e. the normal distribution with mean $5$. Moreover, the weight assigned to ${\rm N}(5,\mbox{0.6})$ differs in the two cases.
The density estimates are displayed in Figure \ref{fig: marginali_peso09_density_100}. The spike corresponding to the common component (concentrated around $5$) is estimated more accurately than the idiosyncratic components (around $0$ and $10$, respectively) of the two samples nicely showcasing the borrowing of information across samples. Moreover, 
the posterior distributions of the number of components are reported in Table \ref{tab:simulated_posterior}. The model correctly detects that each mixture has two components with one of them shared and the corresponding distributions are highly concentrated around the correct values. Finally the Bayes factor $\mathrm{BF}$ to test equality between the two distributions equals $\mbox{0.00022}$ and the null hypothesis of distributional homogeneity is rejected. 

Scenario III consists in generating the data from mixtures with the same components but differing in their weights. Specifically, $\bm{X}_1^{(n_1)}$ and $\bm{X}_2^{(n_2)}$ are drawn from, respectively,
\[
X_{i,1} \sim \mbox{0.8}\, {\rm N}(5,1)+\mbox{0.2}\, {\rm N}(0,1) \quad
X_{j,2} \sim \mbox{0.2}\, {\rm N}(5,1)+\mbox{0.8}\, {\rm N}(0,1),
\]
The posterior distribution of the number of components is again reported in Table~\ref{tab:simulated_posterior} and again the correct number is identified, although in this case the distributions  exhibit a higher variability. The Bayes factor $\mathrm{BF}$ to test equality between the two distributions is 0.54, providing weak evidence in favor of the alternative hypothesis that the distributions differ.

\subsection{Iris dataset}\label{subsec:iris}
Finally, we examine the well known Iris dataset, which
contains several measurements concerning three different species of Iris flower:
setosa, versicolor, virginica. More specifically, we focus on petal width of those species. The sample $\boldsymbol{X}$ has size $n_1=90$, containing $50$ observations of setosa and $40$ of versicolor. The second sample $\boldsymbol{Y}$ is of size $n_2=60$ with $10$ observations of versicolor and $50$ of virginica. 

\begin{figure}[h!]
\begin{center}
\includegraphics[width=0.58\linewidth]{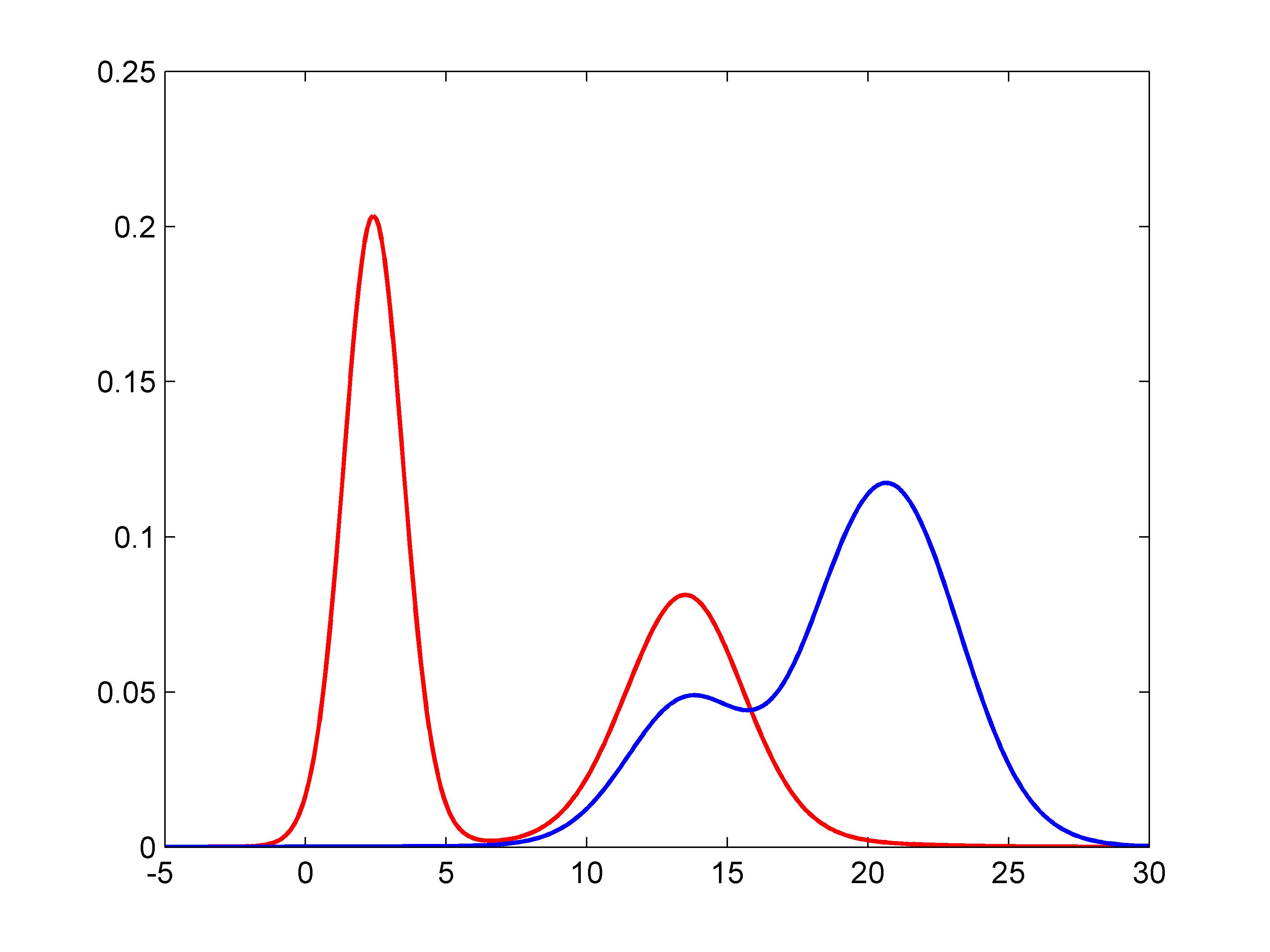}
\end{center}
\begin{flushleft} 
\caption{Estimated densities for $\mathbf{X}$ (red) and $\mathbf{Y}$ (blue).}
\label{fig: iris_densities}
\end{flushleft}
\end{figure}

Since the data are scattered across the whole interval $[0,30]$, we need to allow for large variances and this is obtained by setting $(s_0,S_0)=(1,4)$. The model neatly identifies that the two densities have two components each and that one of them is shared as showcased by the posterior probabilities reported in Table \ref{tab:iris_posterior}.  As for the Bayes factor, we obtain $\mathrm{BF} \approx 0$ leading to the unsurprising conclusion that the two samples
come from two different distributions. The corresponding estimated densities are reported in
Figure \ref{fig: iris_densities}.

\begin{table}[h!]
	\begin{center}
		\begin{tabular}{|l|c|c|c|c|c|c|c|c|}\hline
			 $\#$ comp.  & 0 & 1 & 2 & 3 & 4 & 5 & 6 & $\ge 7$ \\\hline \hline
			$K_1$ 
			& 0 & 0 & \bf{0.466} & 0.307 & 0.141 & 0.055 &  0.020 & 0.011      \\   	
			$K_2$ 
			& 0 & 0.001 & \bf{0.661} & 0.248 & 0.068 & 0.017 & 0.004 & 0.001       \\
			$K_{12}$ 
			& 0  & \bf{0.901} & 0.093  & 0.006 & 0 & 0 & 0 & 0 \\\hline
		\end{tabular}
		\begin{flushleft}
			\caption{Real data: Posterior distributions of the number of components in the first sample
				($K_1$), in the second sample ($K_2$) and shared by the two samples ($K_{12}$). The posterior probabilities corresponding to the MAP estimates are displayed in bold.}
			\label{tab:iris_posterior}
		\end{flushleft}
	\end{center}
\end{table}

\begin{figure}[h!]
	\begin{center}
		\subfigure[]{\includegraphics[width=0.48\linewidth]{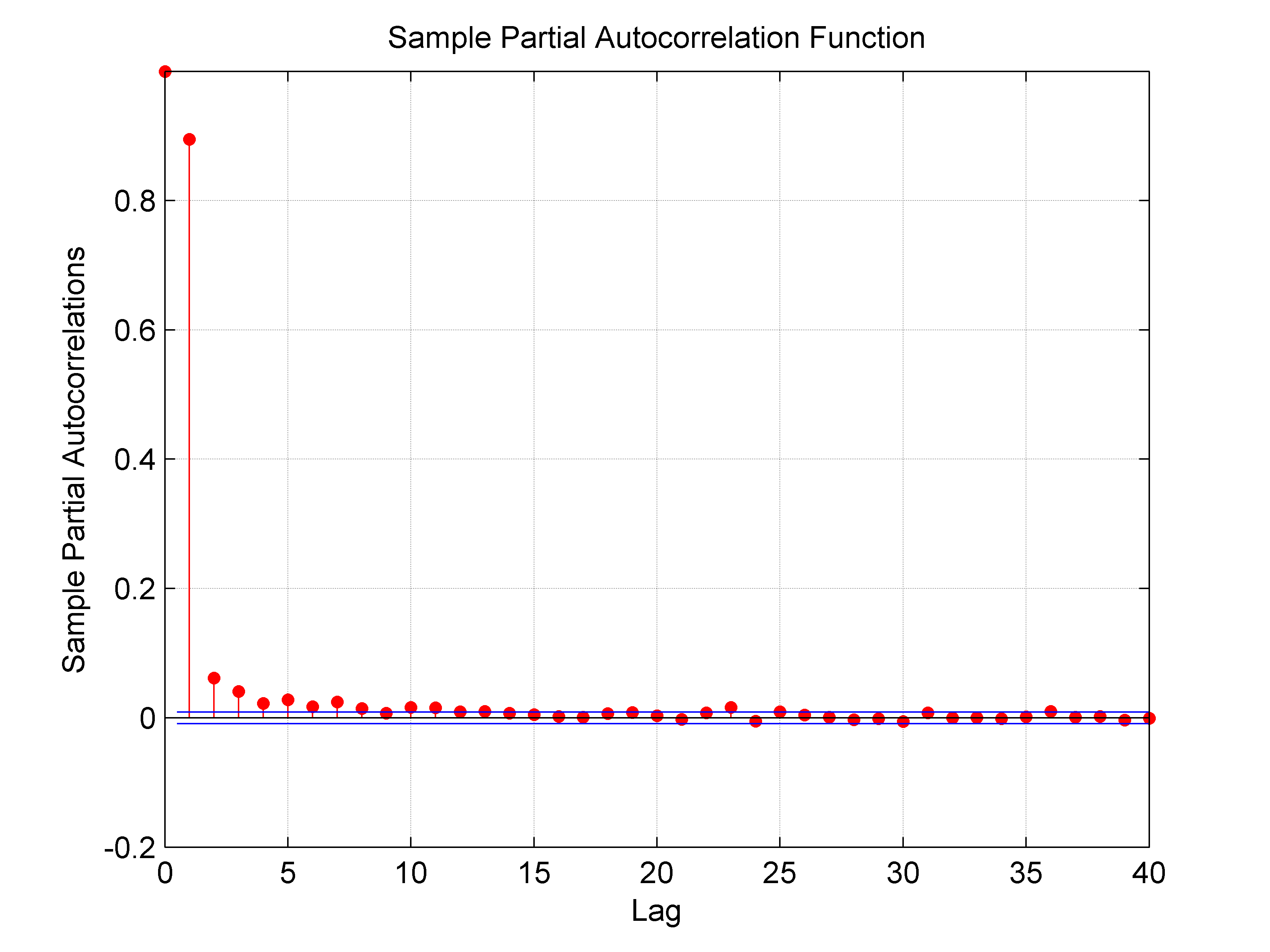}
			\label{fig:parcorr_sigma}}
		\subfigure[]{\includegraphics[width=0.48\linewidth]{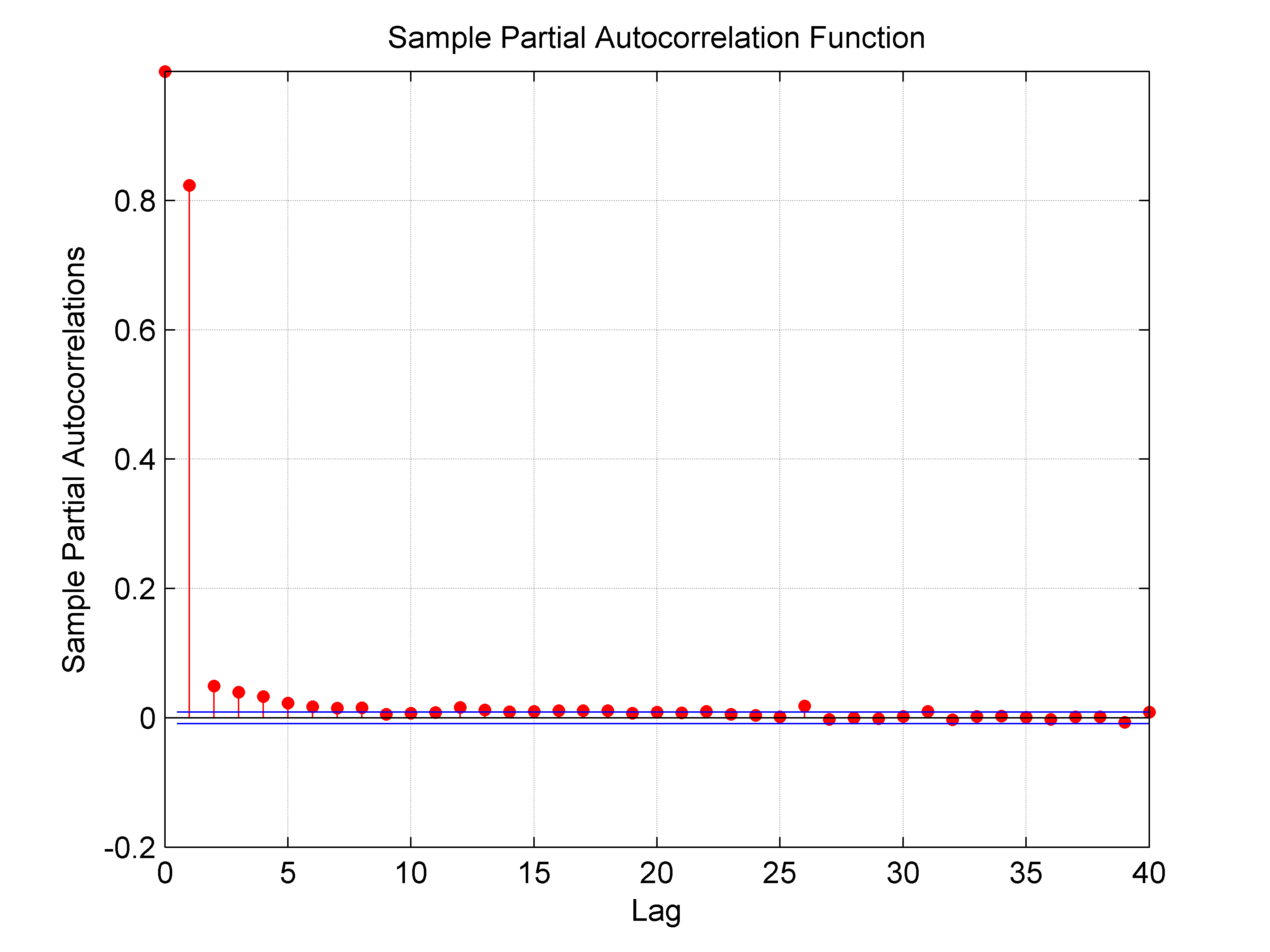}
			\label{fig:parcorr_sigma0}}
	\end{center}
\begin{flushleft}
\caption{Plots of the partial autocorrelation functions for the parameters $\sigma$ (a) and $\sigma_0$ (b).}\label{fig:parcorr}
\end{flushleft}
\end{figure}

We have also monitored the convergence of the algorithm that has been implemented. Though we here provide only details for the Iris dataset, we have conducted similar 
analyses also for each of the illustrations with synthetic datasets in Section~\ref{subsec:synthetic_data}. Notably, all the examples with simulated data have experienced even better performances than those we are going to display henceforth. Figure~\ref{fig:parcorr} depicts the partial autocorrelation function for the sampled parameters $\sigma$ and $\sigma_0$. The partial autocorrelation function apparently has an exponential decay and after the first lag exhibits almost negligible peaks.

\begin{figure}[h!]
	\begin{center}
		\subfigure[]{\includegraphics[width=0.48\linewidth]{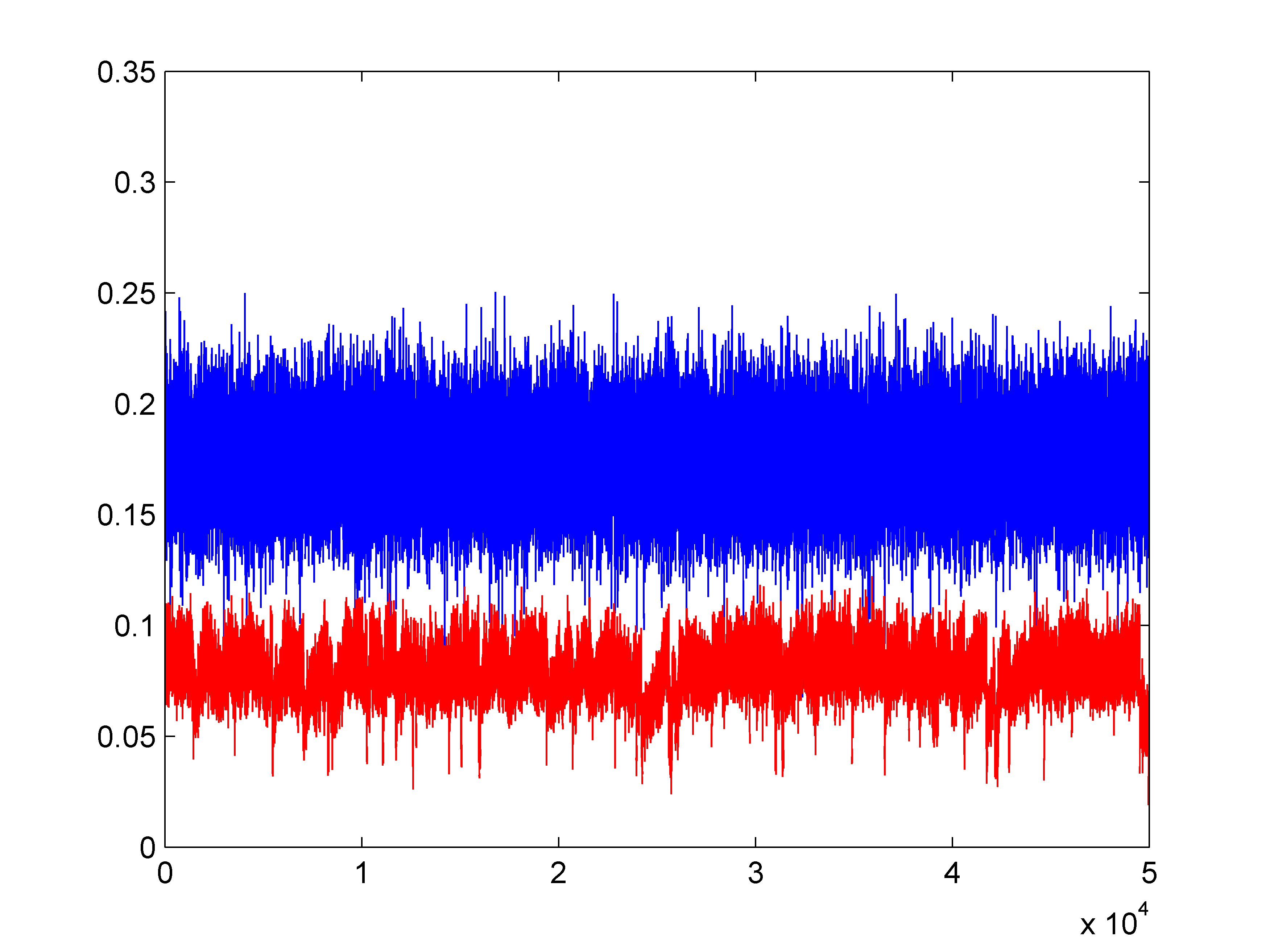}
			\label{fig:convergence_f1}}
		\subfigure[]{\includegraphics[width=0.48\linewidth]{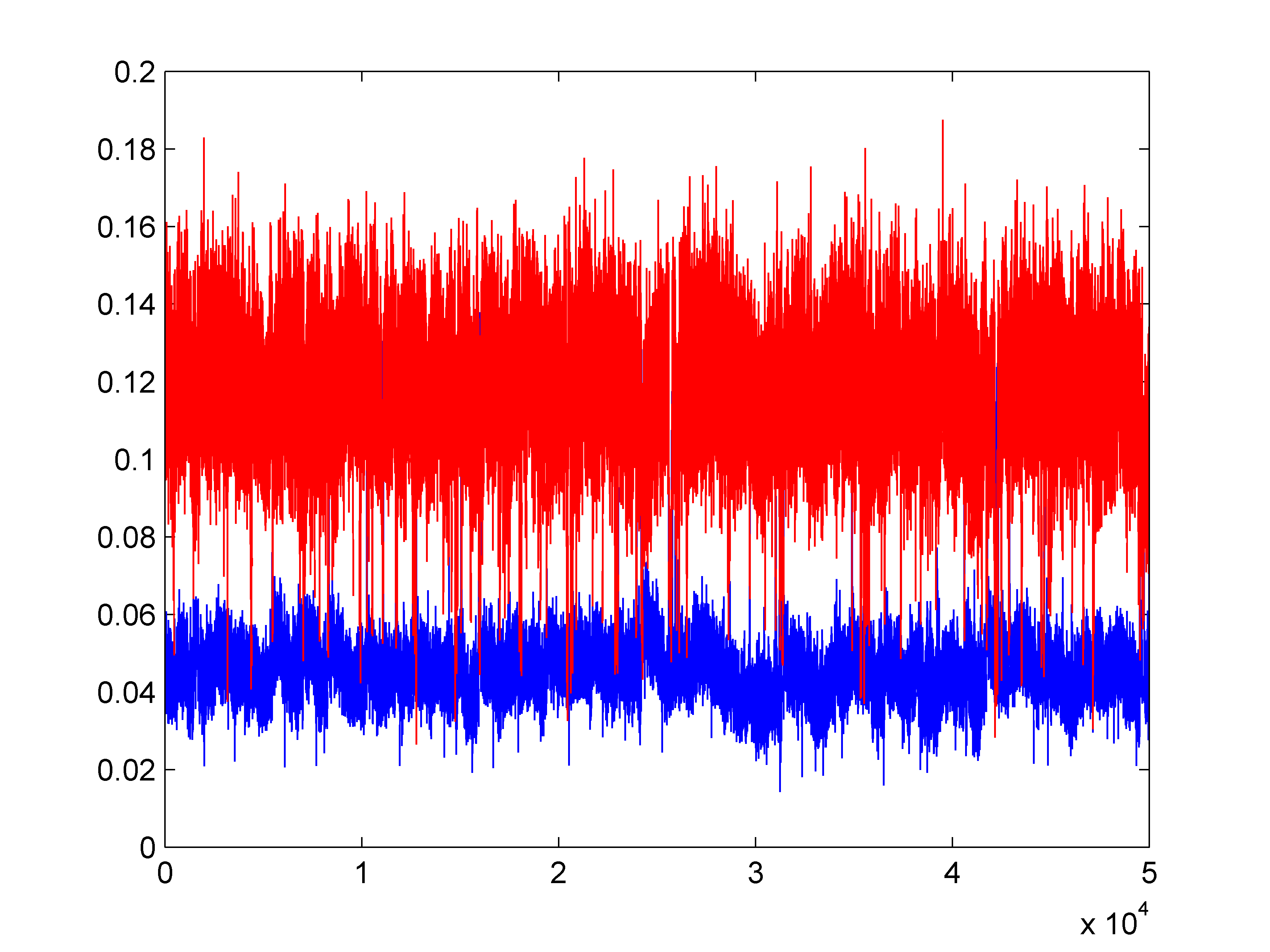}
			\label{fig:convergence_f2}}
	\end{center}
	\begin{flushleft}
		\caption{(a): trace plots of the density referring to $\bm{X}_1^{(n_1)}$ at the points $3$ and $13$; (b): trace plots of the density referring to $\bm{X}_2^{(n_2)}$ at the points $13$ and $21$.}\label{fig:convergence_f}
	\end{flushleft}
\end{figure}

We have additionally monitored the two estimated densities near the peaks, which identify the mixtures' components. More precisely, Figure~\ref{fig:convergence_f1} displays the trace plots of the density referring to the first sample at the points $3$ and $13$, whereas Figure~\ref{fig:convergence_f2} shows the trace plots of the estimated density function of the second sample at the points 
$13$ and $21$.

\section{Concluding remarks}\label{sec:concl}

We have introduced and investigated a novel class of nonparametric priors featuring a latent nested structure. Our proposal allows flexible modeling of heterogeneous data and deals with problems of testing distributional homogeneity in two-sample problems. Even if our treatment has been confined to the case $d=2$, we stress that the results may be formally extended to $d>2$ random probability measures. However, their implementation would be more challenging since the marginalization with respect to $(\tilde p_1,\ldots,\tilde p_d)$ leads to considering all possible partitions of the $d$ random probability measures. While sticking to the same model and framework which has been shown to be effective both from a theoretical and practical point of view in the case $d=2$, a more computationally oriented approach would be 
desirable in this case. There are two possible paths. The first, along the lines of the original proposal of the nested Dirichlet process in \cite{ndp:2008}, consists in using tractable stick--breaking representations of the underlying random probabilities, whenever available to devise an efficient algorithm. The second, which needs an additional significant analytical step, requires the derivation of a posterior characterization of $(\tilde p_1,\ldots,\tilde p_d)$ that allows sampling of the trajectories of latent nested processes and build up algorithms for which marginalization is not needed. Both will be the object of our future research.


\appendix

\section*{Appendix 1}

\subsection*{Proof of Proposition \ref{prp:coincidence}} Since {$(\tilde p_1,\tilde p_2)\sim\mbox{NP}(\nu_0,\nu)$}, one has
\begin{equation}
    \label{eq:porob_coincidence_bis}
    \pi_1
    =\E \:\int_{\P}\tilde{q}^2(\ddr p)=\E
    \:\int_{\P}\frac{\tilde{\mu}^2(\ddr p)}{\tilde{\mu}^2(\P)}=\int_0^\infty u\:\int_{\P}\:\E\,\edr^{-u
        \tilde{\mu}(\P)}\:\tilde{\mu}^2(\ddr p) \:\ddr u
\end{equation}
In order to get the result, we extend and adapt the techniques used in \cite{Jamesall:2006}. Indeed, it can be seen that
\begin{equation}
\E\,\edr^{-u\,\tilde{\mu}(\P_\X)}\:\tilde{\mu}^2(\ddr p)=
\edr^{-c\,\psi(u)}\:\left[c^2\, Q^2(\ddr p)\:\tau_1^2(u)+c\, Q(\ddr p)\,\tau_2(u)\right].
\label{eq:laplace_2}
\end{equation}
Recall that $Q$ is the probability distribution of the NRMI $\tilde{q}_0=\sum_{j\ge 1}\omega_j\:\delta_{\tilde{\theta}_j}$ with $\sum_{j\ge 1}\omega_j=1$ almost surely and $\tilde{\theta}_j\simiid Q_0$. This means that $Q$ is concentrated on the set of discrete probability measures on $\X$. 
If $p=\sum_{j\ge 1}w_j\,\delta_{\theta_j}\in\P$ is fixed, we set $W_{\bm{j},n}:=\{\omega_{j_1}=w_1,\ldots,\omega_{j_n}=w_n\}$ and $\Theta_{\bm{j},n}=\{\tilde{\theta}_{j_1}=\theta_1,\ldots,\tilde{\theta}_{j_n}=\theta_n\}$ where $\bm{j}=(j_1,\ldots,j_n)$ is a vector of positive integer. Then 
\begin{equation}
  Q(\{p\})
  = \P\Big(\tilde q_0=p\Big)\le
  \P\Big[\bigcup_{(*)}\: (W_{\bm{j},n}\cap \Theta_{\bm{j},n})\Big]\label{eq:q_nonatomic}
\end{equation}
where the above union is taken over the set of all vectors $\bm{j}=(j_1,\ldots,j_n)\in\N^n$ such that   $j_1\ne\,\cdots\,\ne j_n$. The upper bound in \eqref{eq:q_nonatomic} is clearly equal to $0$. 
This, combined with \eqref{eq:laplace_2}, yields
\[
\int_{\P_\X}\,\E\,\edr^{-u\,\tilde{\mu}(\P_\X)}\:\tilde{\mu}^2(\ddr p) =c\,\edr^{-c \psi(u)}\,\tau_2(u)\:
\int_{\P_\X}\,Q(\ddr p)=c\,\edr^{-c \psi(u)}\,\tau_2(u)
\]
and the proof is completed. \qed

\subsection*{Proof of Proposition \ref{prp:moment2}}

Let $f_1=\indic_A$ and $f_2=\indic_B$, for some measurable subsets 
$A$ and $B$ of $\P$. One has
\begin{multline*}
 \E \:
 \int_{\P_{\X}^2} f_1(p_1) f_2(p_2) \tilde{q}(\ddr p_1)
\tilde{q}(\ddr p_2)
= \E\: \tilde{q}(A) \tilde{q}(B)\\
=
\E\, \tilde{q}^2(A\cap B) + \E\,\tilde{q}(A \cap B)\tilde{q}(B \cap A^c)+\E\,
\tilde{q}(A \cap B^c)\tilde{q}(B).
\end{multline*}
It can now be easily seen that
\begin{align*}
\E\, \tilde{q}^2(A\cap B)
&=  \E\, \frac{\tilde{\mu}^2(A \cap B)}{\tilde{\mu}^2 (\P_\X)} \\[4pt]
&= \int_0^{\infty}u\, \edr^{-c \psi (u)}\, \Big[c^2 \,Q^2(A \cap B) \,\tau_1^2(u) + c\, Q(A \cap B) \,\tau_2 (u) \Big]\: \ddr u \\[4pt]
&= \pi_1 \,Q^2(A \cap B) +(1-\pi_1)\, Q(A \cap B).
\end{align*}
On the other hand, if $A\cap B=\varnothing$, we get
\[
  \E\, \tilde{q}(A)\tilde{q}(B)= Q(A)Q(B)\: \int_0^{\infty}
  c^2 \, u \,\edr^{-c\, \psi(u)} \:\tau_1^2(u)\, \ddr u =\pi_1 \,Q(A)\,Q(B) 
\]
To sum up, one finds that
\begin{equation*}
\E \,\tilde{q}(A) \tilde{q}(B)
= \pi_1\, Q(A \cap B)+(1-\pi_1)\:[ Q^2(A \cap B)  + Q(A \cap B ) Q(B \cap A^c) 
+Q(A \cap B^c) Q(B)]
\end{equation*}
which boils down to \eqref{eq:2nd_moment}. Now it is easy to prove that \eqref{eq:2nd_moment} is true when $f_1$ and $f_2$ are simple functions and, then, for all positive and measurable functions relying on the monotone convergence theorem. \qed

\subsection*{Proof of Theorem \ref{thm:EPPF_ndp}} 
The partition probability function $\Pi^{(N)}_{k}(
\bm{n}_1,\bm{n}_2,\bm{q}_1,\bm{q}_2)$ equals
\begin{equation}
\label{eq:eppf_integral}
\int_{\X^k} \E  \prod_{j=1}^{k_1} \tilde{p}_1^{n_{j,1}} (\ddr x_{j,1}^*) \prod_{j=1}^{k_2} \tilde{p}_2^{n_{j,2}}
(\ddr x_{j,2}^*) \prod_{j=1}^{k_0} \tilde{p}_1^{q_{j,1}}(\ddr z_j^*) \tilde{p}_2^{q_{j,2}}(\ddr z_j^*) 
\end{equation}
obtained by marginalizing with respect to $(\tilde p_1,\tilde p_2)$. Due to conditional independence of $\tilde{p}_1$ and $\tilde{p}_2$, given $\tilde{q}$, the integrand in \eqref{eq:eppf_integral} can be rewritten as $\E \prod_{\ell=1}^2 h_\ell(\ddr\bm{x}_\ell^*,\ddr \bm{z}^*;\tilde q)$ where, for each $\ell=1,2$, 
\begin{align*}
h_\ell(\ddr\bm{x}_\ell^*,\ddr \bm{z}^*;\tilde q)
&=  \E\Big[ \prod_{j=1}^{k_\ell} \tilde{p}_\ell^{n_{j,1\ell}} (\ddr x_{j,\ell}^*)
\prod_{j=1}^{k_0} \tilde{p}_\ell^{q_{j,\ell}}(\ddr z_j^*) \Big| \:\tilde{q}\: \Big]\\[4pt]
&=\int_{\P_\X} \prod_{j=1}^{k_\ell} p_\ell^{n_{j,\ell}} (\ddr x_{j,\ell}^*)
\prod_{j=1}^{k_0} p_\ell^{q_{j,\ell}}(\ddr z_j^*)\:\tilde q(\ddr p_\ell)
\end{align*}
A simple application of the Fubini--Tonelli theorem, then, yields
\begin{equation}
  \label{eq:eppf_integral_product}
  \Pi^{(N)}_{k}(
\bm{n}_1,\bm{n}_2,\bm{q}_1,\bm{q}_2)=\int_{\X^k} \E \int_{\P_\X^2}\: f_1(p_1)\, f_2(p_2)\:\tilde q(\ddr p_1)\,\tilde q(\ddr p_2)
\end{equation}
where, for each $\ell=1,2$, we have set $f_\ell(p_\ell):=\prod_{j=1}^{k_\ell} p_\ell^{n_{j,\ell}} (\ddr x_{j,\ell}^*)\prod_{j=1}^{k_0} p_\ell^{q_{j,\ell}}(\ddr z_j^*)$
and agree that $\prod_{j=1}^0 a_j \equiv 1$. In view of Proposition~\ref{prp:moment2} the integrand in \eqref{eq:eppf_integral_product} boils down to
\begin{align*}
  \int_{\P_{\X}^2}
  \E  \: f_1(p_1)\, f_2(p_2)\,\tilde q(\ddr p_1)\,\tilde q(\ddr p_2)
  &= \pi_1\,\int_{\P_\X} f_1(p)\, f_2(p)\, Q(\ddr p)+(1-\pi_1)\,\prod_{\ell=1}^2 \int_{\P_\X} f_\ell(p)\, Q(\ddr p)\\[4pt]
  &= \pi_1\,\Big[\E f_1(\tilde {q}_0)\, f_2(\tilde{q}_0)\Big] +(1-\pi_1)\,\Big[\E f_1(\tilde q_0)\Big]\Big[\E f_2(\tilde{q}_0)\Big].
\end{align*}
In order to complete the proof it is now enough to note that, due to non--atomicity of $Q_0$, 
\[
\E f_1(\tilde{q}_0)f_2(\tilde{q}_0)=\E \: \prod_{j=1}^{k_1} \tilde{q}_0^{n_{j,1}} (\ddr x_{j,1}^*)\,
\prod_{j=1}^{k_2} \tilde{q}_0^{n_{j,2}}(\ddr x_{j,2}^*)\,
\prod_{j=1}^{k_0} \tilde{q}_0^{q_{j,1}+q_{j,2}}(\ddr z_j^*)
\]
is absolutely continuous with respect to $Q_0^k$ on $\X^k$ and
\[
\frac{\ddr \E f_1(\tilde{q}_0)f_2(\tilde{q}_0)}{\ddr Q_0^k}(\bm{x}_1^*,\bm{x}_2^*,\bm{z}^*)
=\Phi^{(N)}_{k}(\bm{n}_1,\bm{n}_2,\bm{q}_1+\bm{q}_2)
\]
for any vector $(\bm{x}_1^*,\bm{x}_2^*,\bm{z}^*)$ whose $k$ components are all distinct, and is zero otherwise.
As for the second summand above, from Proposition 3 in \cite{Jamesall:2009} one deduces that
\begin{multline*}
  \Big[\E f_1(\tilde{q}_0)\Big] \,\Big[\E f_2(\tilde{q}_0)\Big]
   = \prod_{j=1}^{k_1} Q_0(\ddr x_{j,1}^*) \, \prod_{j=1}^{k_2} Q_0(\ddr
  x_{j,2}^*) \, \prod_{j=1}^{k_0}Q_0^2(\ddr z_j^*)\\[4pt]
  \times \: \Phi^{(|\bm{n}_1|
    +|\bm{q}_1|)}_{k_0+k_1}(\bm{n}_1,\bm{q}_1) \, \Phi^{(|\bm{n}_2|
    + |\bm{q}_2|)}_{k_0+k_2}(\bm{n}_2,\bm{q}_2) 
\end{multline*}
Then it is apparent that $\left[\E f_1(\tilde{q}_0)\right] \left[\E f_2(\tilde{q}_0)\right]\ll Q_0^k$ and still by virtue of the non--atomicity of $Q_0$ one has
\begin{equation*}
  \frac{\ddr \left[\E f_1(\tilde{q}_0)\right] \left[\E
      f_2(\tilde{q}_0)\right]}{\ddr
    Q_0^k}(\bm{x}_1^*,\bm{x}_2^*,\bm{z}^*)
  =\Phi^{(|\bm{n}_1|
    +|\bm{q}_1|)}_{k_0+k_1}(\bm{n}_1,\bm{q}_1) \, \Phi^{(|\bm{n}_2|
    + |\bm{q}_2|)}_{k_0+k_2}(\bm{n}_2,\bm{q}_2) \:\indic_{\{0\}}(k_0)
\end{equation*}
for any vector $(\bm{x}_1^*,\bm{x}_2^*,\bm{z}^*)\in\X^k$ whose components are all distinct, and is zero otherwise.
Note that if it were $k_0\ge 1$, then some of the infinitesimal factors $Q_0^2(\ddr z_j^*)$ would not cancel and the above density would be exactly equal to zero. \qed

\subsection*{Proof of Proposition \ref{prp:probs_crms}}
Since $\tilde{q}\sim\mbox{NRMI}[\nu;\mathds{M}]$, one has 
$
\tilde{q}=\sum_{j \geq 1}\tilde{\omega}_j \delta_{\tilde{\eta}_j},
$
with $\P (\sum_{j \geq 1} \tilde{\omega}_j  =1 ) = 1$ and $\tilde{\eta}_j \simiid Q$. Furthermore, 
$\tilde{\eta}_j$ is, in turn, a CRM $
\tilde{\eta}_j = \sum_{k \geq 1} \tilde{\omega}_{k}^{(j)} \delta_{X_{k}^{(j)}}
$ 
where
$\P (\sum_{k \geq 1} \tilde{\omega}_{k}^{(j)}  < +\infty )= 1$ and  $X_{k}^{(j)} \simiid Q_0$, for any $j = 1,2 , \dots$. An analogous representation holds true also for 
$\mu_S$, i.e. $\mu_S = \sum_{k \geq 1} \tilde{\omega}_{k}^{(0)} \delta_{X_{k}^{(0)}}$ with the same conditions as above. From the assumptions one deduces that the sequences $( X_{k}^{(j)})_{k \geq 1 }$ and $(\tilde{\omega}_{k}^{(j)})_{k \geq 1}$ are independent also across different values of $j$, and Definition~\ref{def:nlp} entails, with probability 1, 
\[
\P\Big[(\mu_1,\mu_2,\mu_S)\in A_1\times A_2\times A_0\,\Big|\,\tilde q 
\Big]=\tilde q(A_1)\,\tilde q(A_2)\,
\tilde{q}_0(A_0)
\]
which implies  
\begin{align*}
\P(\tilde p_1=\tilde p_2)
  &=
  \E\Big[\P\Big(\frac{\mu_1+\mu_S}{\mu_1(\X)+\mu_S(\X)}=
  \frac{\mu_2+\mu_S}{\mu_2(\X)+\mu_S(\X)}\,\Big|\,\tilde q 
  \Big)\,\Big]\\[6pt]
  &=
  \E \Big[ \sum_{i \not = j} \tilde{\omega}_j \tilde{\omega}_i\: \P \Big(
    \frac{\tilde{\eta}_i +  \mu_S}{ \tilde{\eta}_i (\X) +  \mu_S (\X)} 
    =\frac{\tilde{\eta}_j +  \mu_S}{\tilde{\eta}_j (\X) + \mu_S (\X)} 
    \,\Big|\, \tilde{q} 
    \Big)\\[6pt]
  & \qquad +  \sum_{i=1}^{\infty} \tilde{\omega}_i^2\: \P \Big(\frac{\tilde{\eta}_i +  \mu_S}{\tilde{\eta}_i (\X) +  
  \mu_S (\X)} =\frac{\tilde{\eta}_i + \mu_S}{\tilde{\eta}_i (\X) + \mu_S (\X)} 
    \,\Big|\, \tilde{q} 
    \Big) \Big]  \\[6pt]
  &=\E \Big[ \sum_{i \not = j} \tilde{\omega}_j \tilde{\omega}_i\: \P \Big(
    \frac{\tilde{\eta}_i +  \mu_S}{ \tilde{\eta}_i (\X) + \mu_S (\X)} 
    =\frac{\tilde{\eta}_j + \mu_S}{\tilde{\eta}_j (\X) + \mu_S (\X)} 
    \,\Big|\, \tilde{q} 
    \Big)\Big] +\E \sum_{i=1}^{\infty} \tilde{\omega}_i^2. 
\end{align*}
For the second summand above one trivially has $\E \: \sum_{i=1}^{\infty} \tilde{\omega}_i^2=\P(\mu_1=\mu_2)$. As for the first summand, a simple application of the Fubini--Tonelli theorem and the fact that $\tilde \omega_j\ge 1$, for any $j$, yield the following upper bound 
\begin{align*}
 \E & \Big[\sum_{i \not = j} \tilde{\omega}_j \tilde{\omega}_i\: \P \Big( \frac{\tilde{\eta}_i +  \mu_S}{ \tilde{\eta}_i (\X) + \mu_S (\X)} = 
 \frac{\tilde{\eta}_j +  \mu_S}{\tilde{\eta}_j (\X) +  \mu_S (\X)} 
 \,\Big|\, \tilde{q} 
 \Big)\: \Big] \\[4pt]
 & \quad = 
\sum_{i \not = j} \E \Big[  \tilde{\omega}_j\tilde{\omega}_i \:\P \Big( \frac{\tilde{\eta}_i +  \mu_S}{ \tilde{\eta}_i (\X) + \mu_S (\X)} =
\frac{\tilde{\eta}_j +  \mu_S}{\tilde{\eta}_j (\X) + \mu_S (\X)} 
 \,\Big|\, \tilde{q} 
 \Big)\: \Big] \\[4pt]
  & \quad \leq 
    \sum_{i \not = j}\: \P \Big( \frac{\tilde{\eta}_i +  \mu_S}{ \tilde{\eta}_i (\X) +  \mu_S (\X)} =\frac{\tilde{\eta}_j + \mu_S}{\tilde{\eta}_j (\X) + \mu_S (\X)} \Big).
\end{align*}
The proof is completed by showing that this upper bound is zero. To this end, we fix positive integers $i$, $j$ and $n$, consider the $n$--tuple of atoms $(X_{1}^{(i)}, \cdots , X_{n}^{(i)})$ referring to $\tilde{\eta}_i$ and correspondingly define the sets
\begin{equation*}
\begin{split}
\Theta_{\bm{\ell}}^{(j)} &:= \left\lbrace  \omega \in \Omega : \; X_{1}^{(i)}(\omega) = X_{\ell_1}^{(j)}(\omega), \dots ,
X_{n}^{(i)}(\omega) = X_{\ell_n}^{(j)}(\omega)  \right\rbrace,\\
\Theta_{\bm{\ell}}^{(0)} &:= \left\lbrace  \omega \in \Omega : \; X_{1}^{(i)}(\omega) = X_{\ell_1}^{(0)}(\omega), \dots ,
X_{n}^{(i)}(\omega) = X_{\ell_n}^{(0)} (\omega) \right\rbrace
\end{split}
\end{equation*}
for any $\bm{\ell} = (\ell_1 , \cdots , \ell_n) \in \N^{n}$. It is then apparent that
\[
\P \Big( \frac{\tilde{\eta}_i +  \tilde{\mu}^*_0}{ \tilde{\eta}_i (\X) +  \tilde{\mu}^*_0 (\X)} =\frac{\tilde{\eta}_j +  \tilde{\mu}^*_0}{\tilde{\eta}_j (\X) +  \tilde{\mu}^*_0 (\X)} \Big) 
\leq \P \Big[ \bigcup_{\bm{\ell} \in \N^n , \ell_{h_1} \not = \ell_{h_2} } 
(\Theta_{\bm{\ell}}^{(j)} \cup  \Theta_{\bm{\ell}}^{(0)})
\Big]
\]
and this upper bound is equal to 0, because each of the events $\Theta_{\bm{\ell}}^{(j)}$
and $\Theta_{\bm{\ell}}^{(0)}$ in the above countable union has $0$ probability in view of the 
non--atomicity of $Q_0$ and independence.
\qed

\subsection*{Proof of Theorem \ref{thm:peppf_nested_crms}}

Consider the partition induced by the sample $\bm{X}_1^{(n_1)}$ and $\bm{X}_2^{(n_2)}$ into $k=k_1+k_2+k_0$ groups with frequencies $\bm{n}_\ell=(n_{1,\ell},\ldots,n_{k_\ell,\ell})$, for $\ell=1,2$, and $\bar{\bm{q}}=(q_{1,1}+q_{2,1},\ldots,q_{k_0,1}+q_{k_0,2})$. Recalling that $p_\ell=\mu_\ell/\mu_\ell(\X)$, for $\ell=1,2$, the conditional likelihood is
\begin{equation*}
\prod_{\ell=1}^2\,\prod_{j=1}^{k_\ell}\:p_\ell^{\zeta_{j,\ell}^*n_{j,\ell}}(\ddr x_{j,\ell}^*)\,
{p_S}^{(1-\zeta_{j,\ell}^*)n_{j,\ell}}(\ddr x_{j,\ell}^*)\prod_{r=1}^{k_0}{p_S}^{(1-\zeta_{r,0}^*)(q_{r,1}+q_{r,2})}(\ddr z_r^*)\:
\prod_{\ell=1}^2 p_\ell^{\zeta_{r,0}^*\, q_{r,\ell}}(\ddr z_r^*)
\end{equation*}
where we take $\{x_{j,\ell}:\: j=1,\ldots,k_\ell\}$, for $\ell=1,2$, and $\{z_r^*:\: r=1,\ldots, k_0\}$ as the $k_1+k_2+k_0$ distinct values in $\X$. If we now let
\begin{align*}
  f_0({\mu_S}, u,v) &:=  \edr^{-(u+v){\mu_S}(\X)}\: \prod_{r=1}^{k_0}
{\mu_S}^{(1-\zeta_{r,0}^*)(q_{r,1}+q_{r,2})}(\ddr z_r^*)
\prod_{\ell=1}^2 \,\prod_{j=1}^{k_\ell} \mu_0^{(1-\zeta_{j,\ell}^*)n_{j,\ell}}(\ddr x_{j,\ell}^*)
\\
f_1(\mu_1, u,v) &:= \edr^{-u \mu_1(\X)}\: 
\prod_{j=1}^{k_1} \mu_1^{\zeta_{j,1}^* n_{j,1}} (\ddr x_{j,1}^*)\: 
\prod_{r=1}
^{k_0} \mu_1^{\zeta_{r,0}^* q_{r,1}} (\ddr z_{r}^*) \\
f_2(\mu_2,u,v) &:= \edr^{-v \mu_2 (\X)}\: 
\prod_{j=1}^{k_2} \mu_2^{\zeta_{j,2}^* n_{j,2}} (\ddr x_{j,2}^*)
\prod_{r=1}^{k_0} \mu_2^{\zeta_{r,0}^* q_{r,2}} (\ddr z_{r}^*),
\end{align*}
and further take into account the probability distribution of the labels, conditional on $({\mu_S},\mu_1,\mu_2)$, so that the the joint distribution of the random partition and of the corresponding labels $\bm{\zeta}^{**}=(\bm{\zeta}_1^*,\bm{\zeta}_2^*,\bm{\zeta}_0^*)$ is
\begin{equation} \label{eq:nested_CRM_EPPF2}
\begin{split}
\Pi_k^{(N)}  (\bm{n}_1,\bm{n}_2,\bm{q}_1,\bm{q}_2;\bm{\zeta}^{**})&=
\frac{1}{\Gamma (n_1) \Gamma (n_2)} \int_0^{\infty} \int_0^{\infty} u^{n_1-1} v^{n_2-1} 
\E \Big( \prod_{i=0}^2 f_i(\mu_i, u,v) \Big)\: \ddr u\, \ddr v,
\end{split}
\end{equation}
where, for simplicity, we have set $\mu_0=\mu_S$. 
Now, for any $(u,v)\in\R_+^2$, Proposition~\ref{prp:crm_nested_coincidence} implies 
\begin{align}
\E \: \prod_{i=0}^2 f_i(\mu_i, u,v) 
&= 
\E\: \E \Big(  \prod_{i=0}^2 f_i(\mu_i, u,v) 
\Big| \tilde{q}, \tilde{q}^*\Big) =
\E \: \prod_{i=0}^2 \E \Big(f_i(\mu_i, u,v) \Big| \tilde{q} 
\Big) \nonumber \\
&= 
\Big[\E\: f_0({\mu_S},u,v)\Big]
\E \:\int_{\mathds{M}} f_1(m_1,u,v) f_2(m_2,u,v) \tilde{q}(\ddr m_1) \tilde{q}(\ddr m_2) \nonumber \\
&=\Big[\E\,f_0({\mu_S},u,v)\Big] \label{eq:nested_CRM_EPPF3}
\Big\{ \pi_1^* \Big[\E\, \prod_{i=1}^2 f_i(\tilde{\mu}_0,u,v) \Big]
+(1-\pi_1^*) \prod_{i=1}^2\Big[\E\, f_i(\tilde{\mu}_0,u,v)\Big]  \Big\}
\end{align}
 Using the properties that characterise {$\mu_S$} 
it is easy to show that $\E\: f_0({\mu_S} 
,u,v)\ll Q_0^{k-\bar{k}}$, where $\bar{k}=\sum_{j=1}^{k_1}\zeta_{j,1}^*+\sum_{j=1}^{k_2}\zeta_{j,2}^*+\sum_{r=1}^{k_0}\zeta_{r,0}^*$. Moreover
\begin{equation}
  \label{eq:eppf_labels_zero}
    \frac{\ddr [\E \,f_0({\mu_S} 
    	,u,v)]}{\ddr
      Q_0^{k-\bar{k}}}(\bm{x})= \edr^{-{\gamma}\,
      c_0\,\psi_0(u+v)}\,{\gamma}^{k-\bar{k}}c_0^{k-\bar{k}}\:\prod_{\ell=1}^2\,
    \prod_{j:\,\zeta_{j,\ell}=0}\tau_{n_{j,\ell}}^{(0)}(u+v)\prod_{r=1}^{k_0}\tau_{q_{r,1}+q_{r,2}}(u+v)
\end{equation}
for any $\bm{x}\in\X^{k-\bar{k}}$ with all distinct components, and it is zero otherwise. If one notes that $\prod_{i=1}^2 \E\, f_i(\tilde{\mu}_0,u,v)$ vanishes when at least one of the $\zeta_{r,0}^*$'s is non--zero,
the other terms in \eqref{eq:nested_CRM_EPPF3} can be similarly handled and, after having marginalised with respect to $(\bm{\zeta}_1^*,\bm{\zeta}_2^*,\bm{\zeta}_0^*)$, one has 
\begin{equation*}
\Pi_k^{(N)} (\bm{n}_1 , \bm{n}_2 , \bm{q}_1 , \bm{q}_2) = \pi_1^* \,
\sum_{(\ast\ast )} I_1 (\bm{n}_1 , \bm{n}_2 , \bm{q}_1 + \bm{q}_2, \bm{\zeta}^{**})  
+
(1-\pi_1^*) \sum_{(\ast)} I_2 (\bm{n}_1 , \bm{n}_2 , \bm{q}_1+ \bm{q}_2, \bm{\zeta}^*) 
\end{equation*}
where the first sum runs over all vectors $\bm{\zeta}^{**}=(\bm{\zeta}_1^*,\bm{\zeta}_2^*,\bm{\zeta}_0^*)\in\{0,1\}^k$ and the second sum is over all vectors $\bm{\zeta}^{*}=(\bm{\zeta}_1^*,\bm{\zeta}_2^*)\in\{0,1\}^{k-k_0}$. Moreover,
\begin{align*}
& I_1 (\bm{n}_1 , \bm{n}_2 , \bm{q}_1 + \bm{q}_2, \bm{\zeta}^{**}) = \frac{c_0^k \gamma^{k-\bar{k}}}{
\Gamma(n_1) \Gamma (n_2)} \int_0^{\infty} \int_0^{\infty} u^{n_1-1} v^{n_2 -1} \edr^{-(1+\gamma)c_0 \psi_0(u
+v)}\\ & \qquad\qquad \times \prod_{j=1}^{k_1} \tau_{n_{j,1}}^{(0)} (u+v) 
\prod_{j=1}^{k_2} \tau_{n_{j,2}}^{(0)} (u+v)
\prod_{j=1}^{k_0} \tau_{q_{j,1}+q_{j,2}}^{(0)} (u+v) \ddr u \ddr v. 
\end{align*}
One may further note that
\begin{align*}
\sum_{(\ast )}   I_1 (\bm{n}_1 , \bm{n}_2 , \bm{q}_1 + \bm{q}_2, \bm{\zeta}^{**}) 
&=
\sum_{\bar{k}=0}^{k} \sum_{\left\{ \bm{\zeta}^*  : \; |\bm{\zeta}^* | = \bar{k} \right\}}
\frac{c_0^k \gamma^{k-\bar{k}}}{
	\Gamma(n_1) \Gamma (n_2)} \int_0^{\infty} \int_0^{\infty} u^{n_1-1} v^{n_2 -1} \:
\edr^{-(1+\gamma)c_0 \psi_0(u
	+v)}
\\[4pt] 
& \quad \times \:\prod_{j=1}^{k_1} \tau_{n_{j,1}}^{(0)} (u+v) 
\prod_{j=1}^{k_2} \tau_{n_{j,2}}^{(0)} (u+v)
\prod_{j=1}^{k_0} \tau_{q_{j,1}+q_{j,2}}^{(0)} (u+v) \,\ddr u\, \ddr v \\[4pt]
&= \sum_{\bar{k}=0}^{k} \binom{k}{\bar{k}} \frac{c_0^k \gamma^{k-\bar{k}}}{
	\Gamma(n_1) \Gamma (n_2)}  \int_0^{\infty} \int_0^{\infty} u^{n_1-1} v^{n_2 -1} 
e^{-(1+\gamma)c_0 \psi_0(u
	+v)} \\[4pt] 
& \quad \times\, \prod_{j=1}^{k_1} \tau_{n_{j,1}}^{(0)} (u+v) 
\prod_{j=1}^{k_2} \tau_{n_{j,2}}^{(0)} (u+v)
\prod_{j=1}^{k_0} \tau_{q_{j,1}+q_{j,2}}^{(0)} (u+v) \ddr u \ddr v \\[4pt]
&=  \frac{c_0^k (1+\gamma)^k}{
	\Gamma(n_1) \Gamma (n_2)}  \int_0^{\infty} \int_0^{\infty} u^{n_1-1} v^{n_2 -1} 
e^{-(1+\gamma)c_0 \psi_0(u
	+v)} \\[4pt] 
& \quad \times\: \prod_{j=1}^{k_1} \tau_{n_{j,1}}^{(0)} (u+v) 
\prod_{j=1}^{k_2} \tau_{n_{j,2}}^{(0)} (u+v)
\prod_{j=1}^{k_0} \tau_{q_{j,1}+q_{j,2}}^{(0)} (u+v) \ddr u \ddr v 
\end{align*}
and a simple change of variable yields \eqref{eq:EPPF_nested_CRM2}.
\qed

\subsection*{Details on Examples \ref{exe:stable} and \ref{exe:Dirichlet}}

As for the latent nested $\sigma$--stable process, the first term 
in the expression of the pEPPF \eqref{eq:EPPF_nested_CRM2} turns out to be the EPPF of a 
$\sigma_0$--stable process multiplied by $\pi_1^*=1-\sigma$, namely
\[
(1-\sigma) \frac{\sigma_0^{k-1}\Gamma(k)}{\Gamma(N)}   \prod_{\ell=1}^2 \prod_{j=1}^{k_\ell} 
(1-\sigma_0)_{n_{j,\ell}-1} \prod_{j=1}^{k_0} (1-\sigma_0)_{q_{j,1}+q_{j,2}-1}.
\]
As for the second summand in \eqref{eq:EPPF_nested_CRM2}, the term
$I_2 (\bm{n}_1 , \bm{n}_2 , \bm{q}_1+ \bm{q}_2, \bm{\zeta}^*)$ equals
\begin{multline*}
\frac{\sigma_0^k\gamma^{k-\bar{k}
}}{\Gamma(n_1)
	\Gamma (n_2)} \: \xi_{\sigma_0} (\bm{n}_1 , \bm{n}_2 , \bm{q}_1 + \bm{q}_2)\\[4pt]
\times\:\int_{0}^{\infty} \int_0^{\infty} u^{n_1-1}v^{n_2-1} 
\frac{\exp{ \left\{ -\gamma  (u+v)^{\sigma_0}- u^{\sigma_0}- v^{\sigma_0}\right\}}}{
	(u+v)^{N-\bar{n}_1 -\bar{n}_2 -(k-\bar{k}_1-\bar{k}_2 )\sigma_0} u^{\bar{n}_1 -\bar{k}_1 \sigma_0}
	v^{\bar{n}_2 -\bar{k}_2 \sigma_0}}  \ddr u \ddr v .
\end{multline*}
The change of variables $s=u+v$ and $w=u/(u+v)$, then, yields
\begin{multline*}
I_2 (\bm{n}_1 , \bm{n}_2 , \bm{q}_1+ \bm{q}_2, \bm{\zeta}^*) = \frac{\sigma_0^{k-1}\Gamma(k)\gamma^{k-\bar{k}
}}{\Gamma(n_1)
	\Gamma (n_2)}  \xi_{\sigma_0} (\bm{n}_1 , \bm{n}_2 , \bm{q}_1 + \bm{q}_2)\\
\times\:\int_0^{1} \frac{w^{n_1-\bar{n}_1 +\bar{k}_1 \sigma_0 -1} (1-w)^{
		n_2 -\bar{n}_2 +\bar{k}_2 \sigma_0 -1}}{[ \gamma +w^{\sigma_0}+(1-w)^{\sigma_0}]^k} \ddr w 
\end{multline*}
and the obtained expression for $\Pi_k^{(N)}$ follows.

As far as the latent nested Dirichlet process is concerned, the first term in \eqref{eq:EPPF_nested_CRM2} coincides with the EPPF of a Dirichlet process having total mass $c_0$ multiplied by $\pi_1^*=(c+1)^{-1}$, i.e.
\[
\frac{1}{1+c} \cdot\frac{[c_0 (1+\gamma)]^k}{(c_0(1+\gamma))_N} 
\prod_{\ell=1}^2 \prod_{j=1}^{k_\ell}\Gamma(n_{j,\ell}) \prod_{j=1}^{k_0} \Gamma (q_{j,1}+q_{j,2}).
\]
On the other hand, it can be seen that $I_2 (\bm{n}_1 , \bm{n}_2 , \bm{q}_1+ \bm{q}_2, \bm{\zeta}^*)$ equals
\begin{align*}
&\frac{c_0^k \gamma^{k-\bar{k}
}}{\Gamma(n_1) \Gamma(n_2)}\prod_{\ell=1}^2 \prod_{j=1}^{k_\ell}\Gamma(n_{j,\ell}) \prod_{j=1}^{k_0} \Gamma (q_{j,1}+q_{j,2}) \\
&\qquad\qquad \times
\int_0^{\infty} \int_0^{\infty} \frac{u^{n_1-1}v^{n_2-1}}{(1+u+v)^{\gamma c_0 +N -\bar{n}_1 -
		\bar{n}_2} (1+u)^{\bar{n}_1+c_0} (1+v)^{\bar{n}_2+c_0}} \ddr u \ddr v .
\end{align*}
If ${}_pF_q (\alpha_1, \ldots , \alpha_p ; \beta_1, \ldots , \beta_q; z)$ denotes the generalised hypergeometric series, which is defined as
\[
{}_pF_q (\alpha_1, \ldots , \alpha_p ; \beta_1, \ldots , \beta_q; z) :=
\sum_{k=0}^{\infty} \frac{(\alpha_1)_k \ldots (\alpha_p)_k}{(\beta_1)_k \ldots (\beta_q)_k}
\frac{z^k}{k!},
\]
identity 3.197.1 in \cite{gradshtein1971tables} leads to rewrite $I_2 (\bm{n}_1 , \bm{n}_2 , \bm{q}_1+ \bm{q}_2, \bm{\zeta}^*)$ as follows
\begin{multline*}
c_0^k \gamma^{k-\bar{k}
} \frac{\Gamma((1+\gamma)c_0 +n_1 -\bar{n}_1)}{\Gamma(n_1)\Gamma((1+\gamma)c_0 +N -\bar{n}_1)}\xi_0(\bm{n}_1, \bm{n}_2, \bm{q}_1 +\bm{q}_2) \\[4pt]
\times\:
\int_0^{\infty}  \frac{u^{n_1-1}}{(1+u)^{c_0(1+\gamma)+n_1-\bar{n}_2}}
{}_2F_1(c_0 +\bar{n}_2, n_2; N-\bar{n}_1 +c_0(1+\gamma); -u ) \ddr u .
\end{multline*}
On view of the formula ${}_2F_1 (\alpha, \beta ; \delta; z) = (1-z)^{-\alpha} {}_2F_1 (\alpha , \delta -\beta ; \delta;z/(z-1))$ and of the 
change of variable $t=u/(1+u)$, the integral above may be expressed as  
\begin{equation*}
\int_0^{1}  t^{n_1-1}(1-t)^{c_0(1+\gamma)+c_0-1} 
\times\:
{}_2F_1(c_0 +\bar{n}_2, c_0(1+\gamma) +n_1-\bar{n}_1; N-\bar{n}_1 +c_0(1+\gamma); t ) \ddr t .
\end{equation*}
and, finally, identity 7.512.5 in \cite{gradshtein1971tables} yields the displayed closed form of $\Pi_k^{(N)}$.




\bibliographystyle{agsm}
\bibliography{bib_nested}

\end{document}